%% file: main.tex
\documentclass[a4paper,twoside,10pt]{article}
\usepackage[a4paper,left=2cm,right=2cm, top=3cm, bottom=3cm]{geometry}
\usepackage{cuted}

\usepackage{enumerate}
\usepackage{amsthm}
\usepackage{mathrsfs}
\usepackage{mathtools}
\usepackage{accents}
\usepackage{orcidlink}
\usepackage{graphicx}
\usepackage{amscd,amsmath,amssymb,mathrsfs,bbm,listings}
\usepackage{cancel}
\usepackage{tikz}
\usetikzlibrary{decorations.pathreplacing}
\usepackage{epstopdf}
\graphicspath{{Figs/}}
\usepackage{amsmath}
\usepackage{footmisc}
\usepackage{amsthm}
\usepackage{amssymb}
\usepackage{stmaryrd}
\SetSymbolFont{stmry}{bold}{U}{stmry}{m}{n}
\usepackage{float}
\usepackage{bigints}
\usepackage{cite}
\usepackage{color}
\usepackage[abs]{overpic}
\usepackage[font=footnotesize,labelfont=bf]{caption}
\usepackage{cases}
\usepackage{tikz}
\usepackage{algorithmic}
\usepackage{rotating}
\usepackage{blkarray}
\usetikzlibrary{matrix,calc,arrows}
\usetikzlibrary{pgfplots.groupplots}
\usepackage{pdfcomment}
\usepackage{filecontents}
\usepackage{pgfplots}
\usepackage{ulem}
\usepackage{subcaption}
\usepackage{soul,xcolor}
\usepackage{enumitem}  
\usepackage{verbatim}
\usepackage{graphicx}
\usepackage{alphabeta}
\usepackage{mwe}
\usepackage{algorithm}
\usepackage{array}
\usepackage{adjustbox} 
\newcolumntype{Y}[1]{>{\centering\arraybackslash}m{#1}}

\newcolumntype{C}{>{\centering\arraybackslash}p{2.5cm}}
\numberwithin{equation}{section}

\definecolor{darkblue}{rgb}{0.0, 0.0, 0.8}
\hypersetup{
    colorlinks=true,
    citecolor=red,
    linkcolor=darkblue,
    pdfborder={0 0 0}
}

\newcommand{\tauP}{\texttau{P} }
\newcommand{\Abeta}{A\textbeta{ }}
\mathtoolsset{showonlyrefs,showmanualtags}
\newcommand{\ccrit}{c_h^\mathrm{crit}}

\author{Beatrice Caon\thanks{MOX-Dipartimento di Matematica, Politecnico di Milano, Piazza Leonardo da Vinci 32, 20133, Milano,  Italy (\href{mailto:beatrice.caon@polimi.it}{beatrice.caon@polimi.it})}
\and
Mattia Corti\,\orcidlink{0000-0002-7014-972X}\thanks{MOX-Dipartimento di Matematica, Politecnico di Milano, Piazza Leonardo da Vinci 32, 20133, Milano,  Italy (\href{mailto:mattia.corti@polimi.it}{mattia.corti@polimi.it})}, 
\and
Francesca Bonizzoni\,\orcidlink{0000-0002-6222-3352}\thanks{MOX-Dipartimento di Matematica, Politecnico di Milano, Piazza Leonardo da Vinci 32, 20133, Milano,  Italy (\href{mailto:francesca.bonizzoni@polimi.it}{francesca.bonizzoni@polimi.it})}, 
\and 
Paola F. Antonietti\,\orcidlink{0000-0002-2138-3878}\thanks{MOX-Dipartimento di Matematica, Politecnico di Milano, Piazza Leonardo da Vinci 32, 20133, Milano,  Italy (\href{mailto:paola.antonietti@polimi.it}{paola.antonietti@polimi.it})}, 
}

\title{High-fidelity and Network-based Spatio-temporal Mathematical Models of Alzheimer’s Disease Progression and their Validation Against PET-SUVR Imaging Data\footnote{\textbf{Funding}: This work is partially funded by the European Union (ERC SyG, NEMESIS, project number 101115663). Views and opinions expressed are, however, those of the authors only and do not necessarily reflect those of the European Union or the European Research Council Executive Agency. Neither the European Union nor the granting authority can be held responsible for them. PFA has been partially supported by ICSC—Centro Nazionale di Ricerca in High Performance Computing, Big Data, and Quantum Computing, funded by the European Union—NextGeneration EU. The present research is part of the activities of the Dipartimento di Eccellenza 2023-2027 grant, funded by MUR. 
The authors acknowledge the CINECA award under the ISCRA initiative, for the availability of high-performance computing resources and support under the project IsCc9\_NeuroDG, PI M. Corti, 2025–2026. 
FB is partially funded by  “INdAM - GNCS Project”, codice CUP E53C24001950001.
MC, FB, PFA are members of INdAM-GNCS. }}

\date{}

\begin{document}
\maketitle

\begin{abstract}
Alzheimer's disease is the most common neurodegenerative disorder. Its pathological development is connected with the misfolding and accumulation of two toxic proteins: amyloid-\textbeta{} and tau proteins.
Mathematical models provide a valuable quantitative tool for monitoring disease progression.
In this work, we proposed and compare a novel framework where the spatio-temporal dynamics of amyloid-\textbeta{} and tau proteins is modeled  based on employing either three-dimensional patient-specific geometries or through reduced network-based models defined on the brain connectome.
More specifically, a high-fidelity biophysical model is proposed on three-dimensional brain geometries reconstructed from magnetic resonance imaging, whereas a network-based reduced formulation is defined on the brain connectome. For both approaches, a suitable numerical discretisation is proposed.
A sensitivity analysis is presented to quantify the influence of model parameters on protein concentration patterns as well as compare the quality of the predictions.
For both approaches, the results are validated against PET-SUVR clinical data using [\textsuperscript{18}F]AZD4694 for amyloid-\textbeta{} and [\textsuperscript{18}F]MK6240 for tau protein.
The results indicate that the three-dimensional model provides the most accurate and biologically consistent description of the disease progression, but remains computationally demanding. On the other hand, the reduced graph-based model is cheaper, but it is not always able to achieve reliable results.
\end{abstract}
%
\section{Introduction}
\begin{figure}[t]
    \centering
    \includegraphics[width=0.8\textwidth]{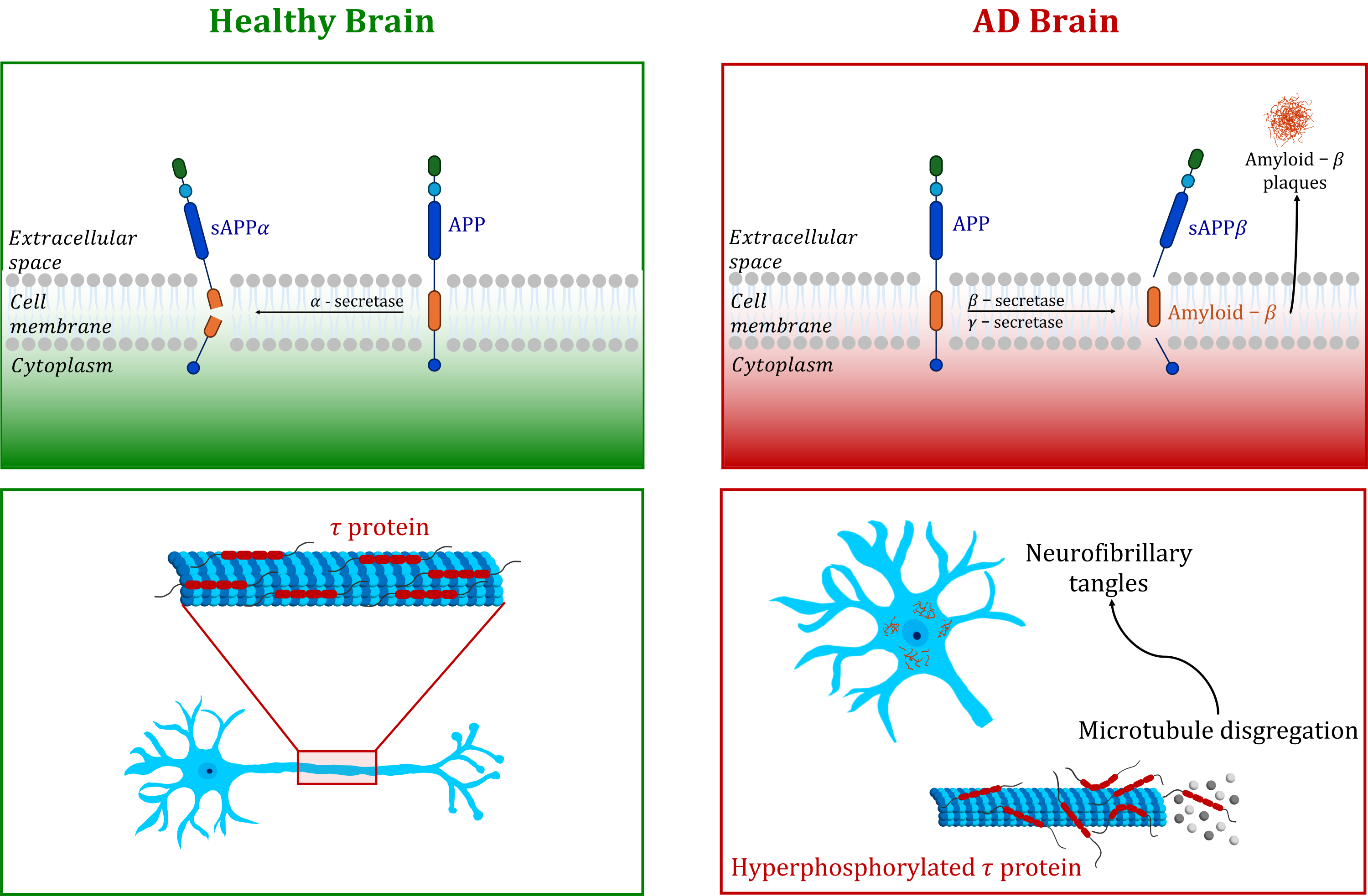}
    \caption{Comparison between healthy and AD brains showing APP processing pathways and \tauP alterations associated with \Abeta and neurofibrillary tangles.}
    \label{fig:Biologia_AD}
\end{figure}

Neurodegenerative diseases represent one of the major medical and societal challenges of our time.
Advances in medicine have significantly reduced mortality at younger ages, contributing to increased life expectancy and shifting the societal burden of disease to older ages \cite{walker_neurodegenerative_2015}.
Some neurodegenerative pathologies are referred to as proteinopathies because they are driven by misfolded proteins that propagate via a \textit{prion-like} mechanism.
Proteins are essential for cellular function and structure. However, when proteins misfold and escape normal degradation processes, they can assemble into pathological aggregates that proliferate and spread within the nervous system \cite{walker_neurodegenerative_2015,jucker_self_2013}. 
The most common proteinopathy is Alzheimer's disease (AD): its pathological development is triggered and sustained by two proteins: amyloid-\textbeta{} (A\textbeta{}) and tau protein (\texttau{P}) \cite{scheltens_alzheimer_2021}. \Abeta is a peptide that derives from the proteolytic processing of amyloid precursor protein (APP), which is present in the brain in physiological conditions. When APP is cleaved by \textalpha{-}secretase, the resulting fragments are non-toxic. In contrast, sequential cleavage by \textbeta{-}secretase followed by \textgamma{-}secretase generates \Abeta peptides prone to misfolding \cite{chen_amyloid_2017} (see Figure~\ref{fig:Biologia_AD}).
Under physiological conditions, \tauP stabilizes microtubules, which are cytoskeletal structures responsible for maintaining neural integrity and supporting intracellular transport. As shown in Figure~\ref{fig:Biologia_AD}, when \tauP becomes hyperphosphorylated, it detaches from microtubules and accumulates in the cytoplasm, where it assembles into filamentous aggregates \cite{jouanne_tau_2017}. An amino terminus in the protein structure facilitates tau-tau interactions and polymer formation \cite{binder_tau_2005}. In AD, the activation of caspase enzymes triggers proteolytic cleavage of \tauP at its carboxy terminus, accelerating the assembly of \tauP into neurofibrillary tangles \cite{binder_tau_2005} (see Figure~\ref{fig:Biologia_AD}). At the earliest stages of the disease, the principal hallmark is the accumulation of \Abeta that induces caspase activation leading to \tauP aggregation.
These misfolded proteins promote neuronal death, leading to impairments in memory, cognition, behavior and/or motor function \cite{wilson_hallmarks_2023}.
The diagnosis of AD remains highly challenging, and in most cases, the pathology is identified only at advanced stages.
For this reason, identifying suitable biomarkers is a crucial objective in research, both for early diagnosis and longitudinal studies \cite{dubois_biomarkers_2023} and for evaluating therapeutic efficacy in clinical trials \cite{pascoal_biomarkers_2024}. 
Biomarkers can be divided into neuroimaging and fluid-based. Among the former, magnetic resonance imaging (MRI) enables \textit{in vivo} assessment of hippocampal and entorhinal cortical atrophy, which are considered characteristic structural alterations of AD \cite{hampel_core_2008}. Positron emission tomography (PET) has gained increasing relevance due to the development of radiotracers that enable \textit{in vivo} monitoring of disease progression by quantifying misfolded protein accumulation. Therefore, PET provides direct insight into the pathological mechanisms underlying AD, although it is an invasive technique due to the use of radiopharmaceuticals, and it is also associated with high costs \cite{bao_pet_2021}.
\par
In this context, mathematical models based on partial differential equations (PDEs) serve as an important tool for the longitudinal monitoring of protein propagation. Among the various modeling approaches one of the simplest yet sufficiently descriptive is the heterodimer model, which focuses exclusively on the interaction between healthy and misfolded protein monomers \cite{matthaus_spread_2009}. Additionally, the heterodimer model can be further reduced to the Fisher-Kolmogorov (FK) equation \cite{fisher_wave_1937,kolmogorov_etude_1937}, widely used to describe biological dynamics. The FK model is particularly suitable under the initial assumption that the concentration of healthy proteins greatly exceeds that of misfolded ones, and its dynamic is negligible \cite{fornari_prion-like_2019}.
Different numerical strategies have been developed to approximate the FK model in the context of neurodegenerative diseases. Starting from the first works based on finite element methods \cite{weickenmeier_physics_2019,schafer_interplay_2019}, recent works have focused on polytopal discontinuous Galerkin (PolyDG) methods \cite{corti_discontinuous_2023,corti_exploring_2023,corti_structure_2024}, also employed in structure-preserving FK formulations to retain qualitative properties of protein propagation dynamics \cite{corti_structure_2024,antonietti_structure_2026, bonizzoni_structure_2020}. Solving such models in three-dimensional geometry is, on the one hand, very informative, but, on the other hand, computationally demanding due to the cost of high spatial resolution.
To overcome this limitation, many studies in the literature either restrict the domain to two-dimensional brain sections \cite{schafer_interplay_2019,corti_exploring_2023} or adopt reduced-order network-based models that represent the brain as a graph derived from the structural connectome \cite{raj_network_2012,iturria-medina_epidemic_2014,fornari_prion-like_2019}, which reflects the preferential diffusion of proteins along axonal pathways. These models rely on several biophysical parameters whose values cannot be directly obtained from clinical data. Hence, parameter estimation has been addressed using inverse techniques and uncertainty quantification methodologies \cite{schafer_bayesian_2021,corti_uncertainty_2024}.\\
\par

\begin{figure}[ht]
    \centering
    \includegraphics[width=\textwidth]{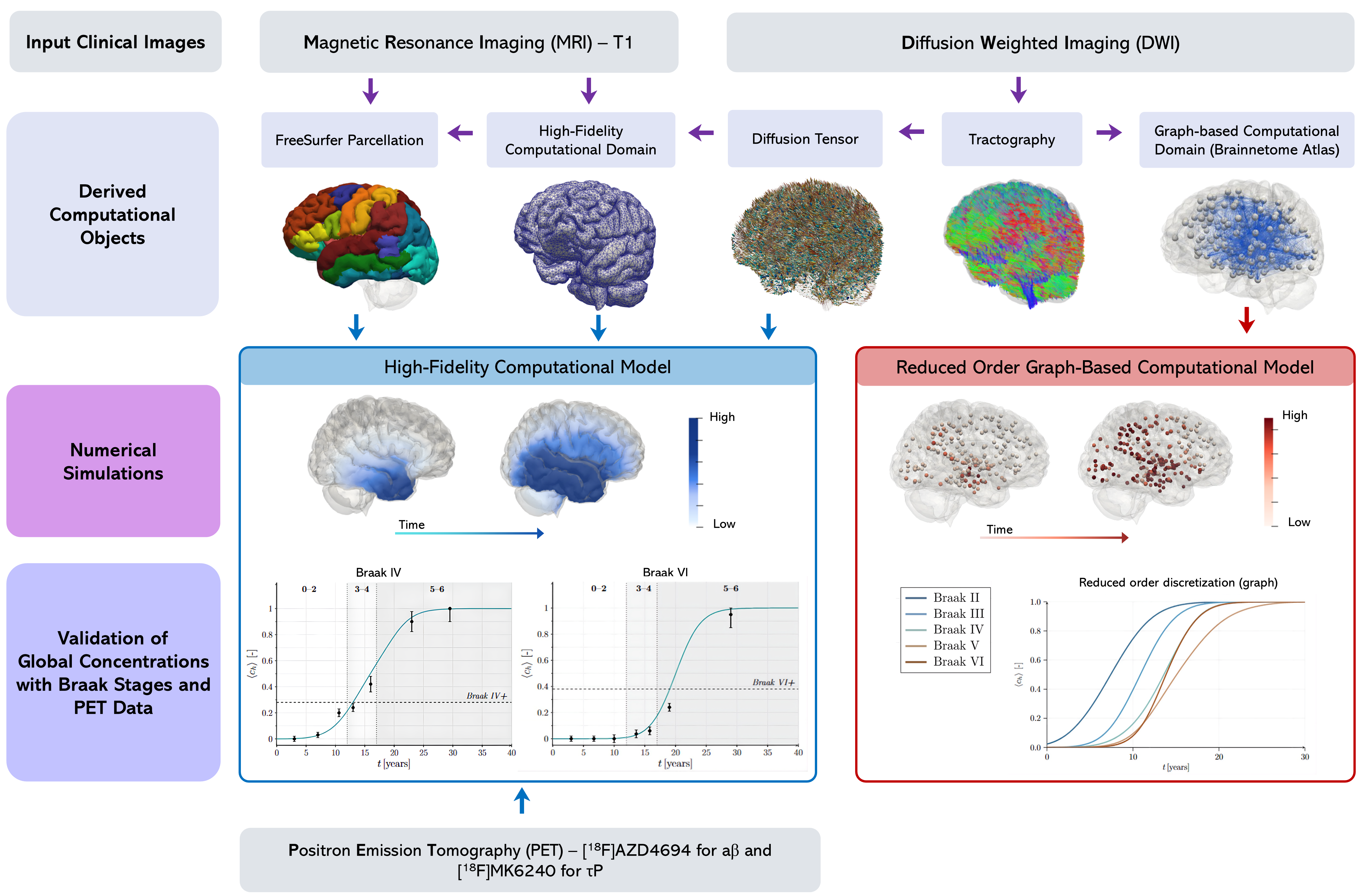}
    \caption{Pipeline from clinical imaging to computational modeling and model validation. Clinical images are used to construct patient-specific brain geometries that define the computational domains for high-fidelity and reduced-order models of protein propagation. Simulated protein concentrations are compared with PET data and Braak staging.}
    \label{BRB_Flowchart}
\end{figure}

In this work, we focus on the FK model to describe the dynamics of \Abeta and \tauP in AD. We compare a high-fidelity discretisation strategy, namely the PolyDG method, defined on three-dimensional brain geometries (3D), with a geometrically reduced two-dimensional sagittal representation (2D) and a reduced-order graph-based formulation (0D) to assess the advantages and limitations of simplified geometric representations of disease progression. Additionally, we analyse the impact of model parameters on protein concentration patterns via a sensitivity analysis, extending the 2D study in \cite{corti_exploring_2023} to the 3D setting and to graph-based models. Finally, we provide a quantitative validation of the 3D model using \textit{in vivo} PET data from \cite{therriault_biomarker_2022}, to evaluate the relevance of physics-based approaches in a clinically oriented context. The measurements are expressed as standardized uptake value ratios (SUVRs) derived from PET imaging with [\textsuperscript{18}F]AZD4694 and [\textsuperscript{18}F]MK6240, which target \Abeta and \tauP, respectively. Attention is also devoted to assessing the ability of the computational models to correctly reproduce the Braak stages of AD pathology \cite{braak_staging_2006}.
Figure~\ref{BRB_Flowchart} summarizes the workflow of the proposed study, from clinical imaging data to computational modeling and validation against PET measurements.\\
\par

The remaining part of the manuscript is structured as follows. In Section~\ref{sec:model}, we introduce the FK mathematical model describing the spatio-temporal spreading of misfolded proteins. Section~\ref{sec:polyDG} presents the high-fidelity \texttt{PolyDG} discretisation method, and Section~\ref{sec:graphs} describes the reduced-order graph-based approach. In Section~\ref{sec:mapping}, we outline the mapping strategy used to validate the numerical simulations against PET data. Section~\ref{sec:results} reports the results of the sensitivity analysis and the validation study. The results and the limitations of the study are discussed in Section~\ref{sec:discussion}. Finally, Section~\ref{sec:conclusions} provides concluding remarks.
\section{A Biophysical Spatio-temporal Mathematical Model: the Fisher-Kolmogorov Equation}
\label{sec:model}
\begin{figure}[htbp]
    \centering
    \begin{subfigure}[t]{0.49\textwidth}
        \centering
        \includegraphics[width=\linewidth]{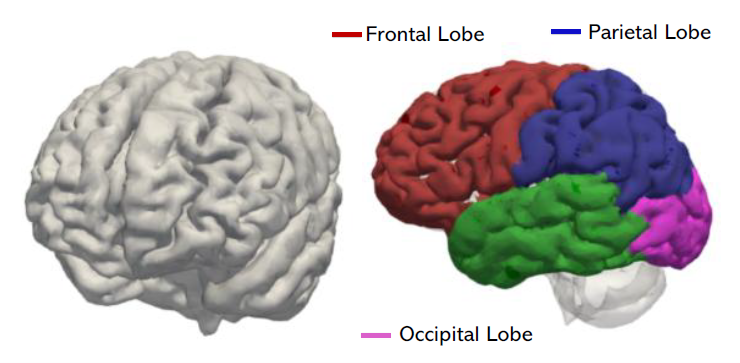}
        \caption{}
        \label{fig:3Dbrain}
    \end{subfigure}
    \hfill
    \begin{subfigure}[t]{0.5\textwidth}
        \centering
        \includegraphics[width=\linewidth]{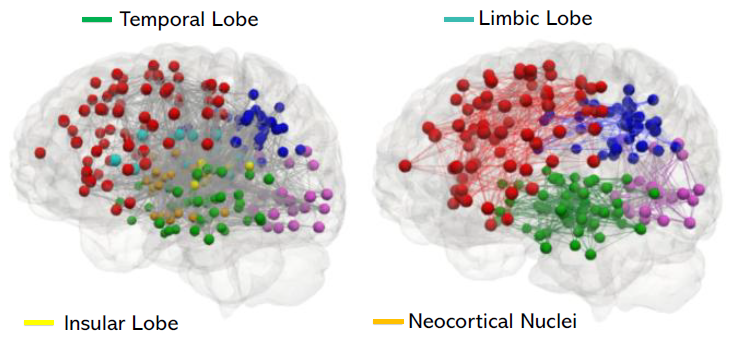}
        \caption{}
        \label{fig:Connectome}
    \end{subfigure}
    
    \caption{Spatial representations of the brain: the high-fidelity brain parenchyma geometry \textbf{(a)} and its reduced connectome-based counterpart \textbf{(b)}. Colors represent different brain regions across both representations.}
    
\end{figure}

In this section, we introduce a biophysical, reaction-diffusion model to describe the production, misfolding, clearance, and spatial spreading of proteins associated with AD. In this work, we adopt the FK equation, as introduced in \cite{weickenmeier_physics_2019,fornari_prion-like_2019}, which can be derived from the heterodimer model \cite{weickenmeier_physics_2019}, exploiting macroscopic considerations of conformational conversion and polymer fragmentation.
The model is formulated at two different spatial resolutions: a physics-based formulation defined over the brain parenchyma (see Figure~\ref{fig:3Dbrain}), and a graph-based formulation defined over the brain connectome (see Figure~\ref{fig:Connectome}).

\subsection{Fisher-Kolmogorov equation on the brain parenchyma}
We consider the brain parenchyma as an open and bounded domain $\Omega\subset \mathbb{R}^d$ with $d\in\{2,3\}$. Given a final time $T > 0$, we seek for the unknown relative concentration of misfolded protein $ c = c(\boldsymbol{x}, t)\in [0,1]$ such that
\begin{equation}
\label{eq:fisher_kolmogorov}
\frac{\partial c}{\partial t}
= \nabla \cdot (\mathbf{D} \cdot \nabla c)
+ \alpha c \left(1 - c \right)
\quad \textrm{in}\ \Omega\times(0,T].
\end{equation}
Here, $c\equiv0$ corresponds to the absence of misfolded proteins, whereas $c\equiv1$ denotes a high prevalence of them. 
In equation \eqref{eq:fisher_kolmogorov}, the diffusion tensor $\mathbf{D} = d_{\mathrm{ext}} \mathbf{I} + d_{\mathrm{axn}} \, \mathbf{\bar{a}} \otimes \mathbf{\bar{a}}$
characterizes the protein spreading \cite{weickenmeier_physics_2019}.
Specifically, $d_\mathrm{ext}$ models the extracellular diffusion and the term $d_\mathrm{axn}$ models the anisotropic diffusion along the axonal direction, denoted by the unit vector $\bar{\mathbf{a}} = \bar{\mathbf{a}}(\boldsymbol{x})$. These directions can be derived from diffusion weighted images (DWI) \cite{corti_discontinuous_2023}. The parameter $\alpha$ denotes the conversion rate of proteins from healthy to misfolded, and it represents the biological processes of production, misfolding, and clearance.
To complete the mathematical model, we supplement the PDE~\eqref{eq:fisher_kolmogorov} with suitable initial and boundary conditions, leading to the following boundary value problem in strong form:
for all $t \in (0,T]$, find $c=c(\boldsymbol{x},t)\in[0,1]$ such that
\begin{gather}
\label{eq:strong_fkpp}
\left\{
\begin{aligned}
    &\frac{\partial c}{\partial t} = \,\nabla \cdot (\mathbf{D} \nabla c) + \alpha c (1-c),  & \quad \mathrm{in}\ \Omega \times (0,T], \\[1ex]
    &(\mathbf{D} \nabla c) \cdot \boldsymbol{n}_\Omega = \, 0, & \quad \mathrm{on}\ \partial \Omega \times (0, T], \\[1ex]
    &c(\cdot,0) = \, c_0, &\quad\mathrm{in}\ \Omega,
\end{aligned}
\right.
\end{gather}
where $\boldsymbol{n}_\Omega$ denotes the outward normal at the boundary $\partial \Omega$ and $c_0=c_0(\boldsymbol{x})$ denotes the initial condition.
Due to the unstable nature of the equilibrium $c = 0$ and the stability of the equilibrium $c = 1$, any positive initial concentration $c_0$ leads to progressive protein accumulation. The homogeneous Neumann boundary condition $(\mathbf{D} \nabla c) \cdot \boldsymbol{n}_\Omega = \, 0$ reflects the assumption that the brain parenchyma does not exchange misfolded proteins with the surrounding cerebrospinal fluid.

\subsection{Fisher-Kolmogorov equation on the brain connectome}
\label{sec:graph_formulation}
In AD, misfolded proteins propagate from localized infection regions along axonal fiber tracts throughout the entire brain \cite{weickenmeier_physics_2019}.
Since the spreading is constrained by the structural connectivity, the FK equation \eqref{eq:fisher_kolmogorov} can be projected onto the brain connectome, reducing the spatial complexity.
We model the brain connectome as a connected network $\Gamma = \bigcup_{k=1}^{N_\Gamma} \gamma_k$, where each $\gamma_k$ represents the $k$-th axonal tract and $N_\Gamma$ denotes the total number of tracts. Fibers intersect only at their endpoints, which define the network nodes.
The FK equation on the brain connectome reads: for all $t \in (0,T]$, find the relative concentration of misfolded protein along the $k$-th axonal tract $c_k(s,t)\in [0,1]$ such that
\begin{gather}
\label{eq:graph_fkpp}
\left\{
\begin{aligned}
    &\frac{\partial c_k(s,t)}{\partial t} 
= \frac{\partial}{\partial s} \Big( D_k \frac{\partial c_k(s,t)}{\partial s} \Big) + \alpha c_k(s,t) (1 - c_k(s,t)), & \mathrm{in} \ \gamma_k \times (0,T], \ k = 1,\dots N_\Gamma , \\[1em]
&c_k(s,0) = c_{0,k}(s) & \mathrm{in}  \  \gamma_k,\ k = 1,\dots N_\Gamma,
\end{aligned}
\right.
\end{gather}
where $s$ is the curvilinear coordinate along the axon, and $D_k$ is the effective diffusivity along the $k$-th axonal tract.
At each graph node, we assume continuity of the protein concentration and the conservation of the flux across connecting fibers.
\section{Numerical discretisation }
%
In this section, we first introduce a high-fidelity discretisation method for problem~\eqref{eq:strong_fkpp} (see Section~\ref{sec:polyDG}), and then we discuss a graph-based discretisation method for problem~\eqref{eq:graph_fkpp} (see Section~\ref{sec:graphs}).
\subsection{High-fidelity discretisation: PolyDG method}
\label{sec:polyDG}
Here, we recall the PolyDG discretisation method introduced in \cite{corti_discontinuous_2023}. First, we introduce a polytopal mesh partition $\mathscr{T}_h$ of the domain $\Omega$ made of disjoint polygonal/polyhedral elements $K$, and we denote with $\mathscr{F}_h^\mathrm{I}$ the set of ($d-1$)--dimensional internal faces.
We define $\mathbb{P}_\ell(K)$ as the space of polynomials of degree $\ell \geq 1$ over the mesh element $K$ and $W_{h}^\mathrm{DG} = \left\{\, 
w \in L^2(\Omega) : \; w|_{K} \in \mathbb{P}_{\ell}(K) 
\;\; \forall K \in \mathscr{T}_h 
\right\}$ as the discontinuous finite element space.
For each internal face $F \in \mathscr{F}_h^\mathrm{I}$ shared by elements $K^\pm$, we denote by $\mathbf{n}^\pm$ the outward unit normals of $K^\pm$ on $F$. For sufficiently regular scalar-valued function $v$ and a vector-valued function $\mathbf{q}$, the trace operators are defined as \cite{arnold_unified_2002}:
\begin{alignat}{3}
\{\!\!\{v\}\!\!\} &= \frac{1}{2}(v^+ + v^-), &\quad \llbracket v \rrbracket &= v^+ \mathbf{n}^+ + v^- \mathbf{n}^-, &\quad \text{on } F \in \mathscr{F}_h^\mathrm{I}, \\
\{\!\!\{\mathbf{q}\}\!\!\} &= \frac{1}{2}(\mathbf{q}^+ + \mathbf{q}^-), &\quad \llbracket \mathbf{q} \rrbracket &= \mathbf{q}^+ \cdot \mathbf{n}^+ + \mathbf{q}^- \cdot \mathbf{n}^-, &\quad \text{on } F \in \mathscr{F}_h^\mathrm{I}.
\end{alignat}
Here, the superscripts $\pm$ indicate the traces of functions on $F$ taken from the interiors of $K^\pm$, respectively.
\par
We define the bilinear form $\mathscr{A}_h :  W_h^\mathrm{DG} \times W_h^\mathrm{DG} \to \mathbb{R}$ as:
\begin{equation*}
\mathscr{A}_h(u_h,v_h) = 
\int_{\Omega} (\mathbf{D} \nabla_h u_h) \cdot \nabla_h v_h 
+ \int_{\mathscr{F}_h^\mathrm{I}} 
\Big( \eta \llbracket u_h \rrbracket \cdot \llbracket v_h \rrbracket 
- \{\!\!\{ \mathbf{D} \nabla_h u_h \}\!\!\} \cdot \llbracket v_h \rrbracket 
- \{\!\!\{ \mathbf{D} \nabla_h v_h \}\!\!\} \cdot \llbracket u_h \rrbracket \Big), \quad \forall u_h,v_h \in W_h^\mathrm{DG},
\end{equation*}
where $\nabla_h$ denotes the elementwise gradient and $\eta: \mathscr{F}_h^\mathrm{I} \to \mathbb{R}_+$ is the penalization function defined as in \cite{corti_discontinuous_2023}:
\begin{equation}
    \eta = \frac{\ell^2}{\{ h \}_H}\eta_0 
    \max \left( \{ d^K \}_H, \alpha \right),
\end{equation}
with $\eta_0$ being a constant, $d^K = \bigl\| \sqrt{\mathbf{D}|_K} \bigr\|^2$ and $\{ \cdot \}_H$ denotes the harmonic average operator,
defined by $\{ v \}_H = \frac{2 v^+ v^-}{v^+ + v^-}$, and $\|\cdot\|$ is the supremum norm.
The semi-discrete PolyDG formulation of~\eqref{eq:strong_fkpp} reads: find $c_h=c_h(t) \in W_h^\mathrm{DG}$ such that, for all $t \in (0,T]$, there holds
\begin{gather}
\label{eq:weak_fkpp}
\left\{
\begin{aligned}
  &\left( \partial_t c_h, w_h \right)_{\Omega} 
+ \mathscr{A}_h(c_h, w_h) - \left(\alpha c_h, w_h \right)_\Omega + \left(\alpha c_h^2, w_h \right)_\Omega = \, 0, && \quad \forall\, w_h \in W_h^\mathrm{DG}, \\
& c_h(\cdot,0) = \, c_h^0, && \quad\mathrm{in}\ \Omega,  
\end{aligned}
\right.
\end{gather}
where $c_h^0 \in W_h^\mathrm{DG}$ is a suitable PolyDG approximation of the initial condition $c_0$. 
To derive the algebraic form of~\eqref{eq:weak_fkpp}
we introduce a basis $\{\varphi_j\}_{j=1}^{N}$ of $W_h^\mathrm{DG}$, with $N=\textrm{dim}(W_h^\mathrm{DG})$, and we express the solution as linear combination of basis functions $\{\varphi_j\}_{j=1}^{N}$ with time dependent coefficients $\{C_j\}_{j=1}^{N}$:
\begin{equation}
\label{eq:sol}
   c_h(\boldsymbol{x},t) = \sum_{j=0}^{N} C_j(t) \varphi_j(\boldsymbol{x}).
\end{equation}
Choosing the test function $w_h$ in~\eqref{eq:weak_fkpp} as a basis function $\varphi_i$ and using~\eqref{eq:sol}, we derive: 
\begin{subequations}
\label{eq:algebraic_fkpp}
\begin{alignat}{3}
\mathbf{M}\dot{C}(t) + \mathbf{A}C(t) - \mathbf{M}_\alpha C(t) + \widetilde{\mathbf{M}}_\alpha(C(t))C(t) = & \,\boldsymbol{0}, 
&& \qquad t \in (0,T], \\
C(0) = &\, C_0,
\end{alignat}
\end{subequations}
where $\mathbf{M},\, \mathbf{A},\, \mathbf{M}_\alpha$ and $\mathbf{\widetilde{M}}_\alpha$ are square matrices of size $N$ defined as:
\begin{gather}
\label{eq:matrices}
\begin{aligned}
[\mathbf{M}]_{ij} = & \, (\varphi_j, \varphi_i)_{\Omega} && \qquad \text{(Mass matrix)} \\
[\mathbf{A}]_{ij} = & \, \mathcal{A}_h(\varphi_j, \varphi_i) && \qquad \text{(Stiffness matrix)}, \\
[\mathbf{M}_\alpha]_{ij} = & \, (\alpha \varphi_j, \varphi_i)_{\Omega} && \qquad \text{(Linear reaction matrix)} \\
[\mathbf{\widetilde{M}}_\alpha(C(t))]_{ij} = & \, (\alpha c_h(t)\varphi_j, \varphi_i)_{\Omega} 
&& \qquad \text{(Nonlinear reaction matrix)}
\end{aligned}
\end{gather}
and $C_0$ is the vector of length $N$ representing  $c_h^0$.
\par

To discretize in time, we introduce a uniform partition of $[0,T]$ in $N_T$ intervals $I_n=(t_{n-1},t_n]$ for $n=1,\ldots,N_T$ with $t_0=0$ and $t_n=t_0+n\Delta t$.
Applying the semi-implicit Euler method to~\eqref{eq:algebraic_fkpp}, we get the following fully discrete formulation of~\eqref{eq:strong_fkpp}: 
for $n=0,\dots,N_T-1$, find $C^{n+1} \simeq C(t_{n+1})$ such that
\begin{equation}
\left(\mathbf{M} + \Delta t\,\left(\mathbf{A} - \mathbf{M}_\alpha + \mathbf{\widetilde{M}_\alpha}(C^n)\right)\right)C^{n+1}
= \mathbf{M}C^n,
\end{equation}
where $C^0=C_0\in\mathbb{R}^N$. 
A detailed discussion on the method presented here can be found in \cite{corti_discontinuous_2023}.
\subsection{Reduced-order discretisation: graph-based finite differences}
\label{sec:graphs}
In this section, we propose a discretisation method for the connectome-based formulation~\eqref{eq:graph_fkpp} \cite{betzel_generative_2017,corti_uncertainty_2024,fornari_prion-like_2019}. 
For computational purposes, the brain connectome $\Gamma = \bigcup_{k=1}^{N_\Gamma} \gamma_k$ can be equivalently represented as a weighted graph $\mathscr{G} = (V,E)$, where $E$ is the set of edges corresponding to the axonal tracts $\{\gamma_k\}_{k=1}^{N_\Gamma}$ and $V = \{\mathbf{x}_j\}_{j=1}^M$ is the set of nodes being the endpoints of $\{\gamma_k\}_{k=1}^{N_\Gamma}$.
Each edge connecting two nodes $i$ and $j$ is associated with a weight $w_{ij}$. 
Adopting the PDEs discretisation framework on graph \cite{smola_kernels_2003}, we obtain the following semi-discrete formulation: find $\{C_j(t)\}_{j=1}^M$ such that, for all $t \in (0,T]$ there holds
\begin{gather}
\label{eq:graph_fkpp_discretisation}
\left\{
\begin{array}{ll}
    \dot{C}_j(t) &= - \sum_{i=1}^{M} \mathbf{L}_{ji} \, C_i(t) + \alpha\, C_j(t)\,(1 - C_j(t)),\\
    C_j(0) &= C_{0j},
\end{array}
\right.
\end{gather}
where $C_j(t)$ is the approximation of the solution in the $j$-th node $\mathbf{x}_j$ at time $t$, namely $C_j(t)\simeq c(\mathbf{x}_j,t)$, and $C_{0j}$ is the approximation of the initial condition in $\mathbf{x}_j$, namely $C_{0j}\simeq c_0(\mathbf{x}_j)$. 
The matrix $\mathbf{L} \in\mathbb{R}^{M\times M}$ is the Laplacian matrix, representing the discrete analogue of the continuous Laplacian operator \cite{smola_kernels_2003}. It is expressed as $\mathbf{L} = \mathbf{K} - \mathbf{A}$, where $\mathbf{K}$ is the degree matrix - a diagonal matrix measuring the strength of each node - and $\mathbf{A}$ is the adjacency matrix of the graph weights $\{w_{i,j}\}_{i,j=1}^M$. More in details, the elements of the Laplacian matrix are:
\begin{equation*}
 \mathbf{L}_{ij} =
\begin{cases}
    -w_{ij}, & i \neq j, \\[6pt]
    \displaystyle\sum_{k=1}^{M} w_{ik}, & i = j.
\end{cases}   
\end{equation*}
\par
To discretize in time, we apply a time stepping method to~\eqref{eq:graph_fkpp_discretisation}. In particular, we employ the Crank-Nicolson scheme with a semi-implicit treatment of the nonlinear reaction term using a second-order extrapolation as in \cite{corti_uncertainty_2024}.  Let $0 = t_0 < t_1 < \cdots < t_{N_T} = T$ be a uniform partition of $[0,T]$ with time step size $\Delta t = T/N_T$. The fully discrete formulation of problem~\eqref{eq:graph_fkpp} reads: given the initial conditions $C^0 = C_0$ and $C^{-1} = C_{-1}$, find $\{C^{n+1}\}_{n=0}^{N_T-1}$, such that:
\begin{subequations}
\label{eq:graph__discrete_fkpp}
\begin{alignat}{3}
\nonumber 2C^{n+1} + & \, \Delta t\mathbf{L} C^{n+1} - \alpha \Delta t \left( \mathbf{1} - \frac{1}{2} \left(3C^n - C^{n-1}\right)\right) \odot C^{n+1}
\\ = & \, 2 C^n - \Delta t \mathbf{L} C^n - \alpha \Delta t \left( \mathbf{1} - \frac{1}{2} \left(3C^n - C^{n-1}\right)\right) \odot C^n,
& & 
\end{alignat}
\end{subequations}
where $\mathbf{1}$ is a vector of ones, and the symbol $\odot$ denotes the componentwise Hadamard product. In the adopted notation, $C^{n+1}$ is a vector of $M$ components $\{C^{n+1}_j\}_{j=1}^M$, and its $j$-th element approximates the solution at node $\mathbf{x}_j$ and time instance $t_n$.
\subsection{Simulation setup}
In this section, we describe the simulation setup for both the high-fidelity and graph-based computations. 
\subsubsection*{High-fidelity simulations}
For the numerical simulations, we construct a 3D brain mesh from a structural MRI taken from the OASIS-3 database \cite{lamontagne_oasis_2019}. The structural MRI is segmented using FreeSurfer \cite{fischl_freesurfer_2012} with the Desikan--Killiany atlas \cite{desikan_automated_2006}. The resulting surfaces are processed with the SVMTK library to generate a volumetric unstructured mesh composed of $142\,658$ elements \cite{mardal_mathematical_2022}. 
The axonal directions $\overline{\mathbf{a}}$ to construct the diffusion tensor $\mathbf{D}$ are obtained from DWI data using FreeSurfer and the Nibabel library. Since the protein spreading predominantly occurs along the axonal directions, the axonal diffusion is set to be ten times faster than the isotropic diffusion, with $d_{\mathrm{ext}} = 8\ \mathrm{mm}^2/\mathrm{year}$ and $d_{\mathrm{axn}} = 80\ \mathrm{mm}^2/\mathrm{year}$ \cite{corti_discontinuous_2023}.
Regarding the discretisation parameters, we fix the polynomial order $\ell = 2$ in each mesh element and set the penalty parameter $\eta_0 = 10$. Moreover, we adopt a uniform time step $\Delta t = 0.05$ years, and we simulate up to a final time $T=40$ years. Concerning the initial conditions, \Abeta and \tauP are initially located in the cerebral cortex and in the enthorinal cortex, respectively. Simulations were run on the GALILEO100 supercomputer (528 nodes, each with $2\times$Intel Cascade Lake 8260 CPUs, 24 cores, 2.4~GHz, 384~GB RAM) at the CINECA supercomputing centre.
\subsubsection*{Graph-based simulations}
The brain connectivity graph is derived from DWI data using DSI Studio tractography, which reconstructs the principal orientations of axonal fibres in the brain \cite{yeh_dsi_2025}. A weighted graph is then constructed, where nodes correspond to brain regions defined by the Brainnetome atlas \cite{fan_human_2016}. The weights of the graph edges are computed as $w_{ij} = k {n_{ij}}/{l_{ij}}$ for $i,j = 1,\dots,M$, where $n_{ij}$ denotes the number of tracts that connect two brain regions associated with different nodes $i$ and $j$, while $l_{ij}$ is the mean value length of these tracts. The coefficient $k$ is a scaling factor calibrated to obtain consistent times with previous literature works \cite{fornari_prion-like_2019,schafer_network_2020,corti_uncertainty_2024}.
In all numerical simulations, we adopt a uniform time step $\Delta t = 0.05\ \mathrm{year}$ and simulate up to a final time $T = 40\ \mathrm{years}$. Initial conditions are defined to reflect the characteristics of the pathology, namely, \Abeta proteins are initially imposed in the cerebral cortex and \tauP seeding is assigned to nodes corresponding to the entorhinal cortex and the rostral portion of the parahippocampal gyrus.
\section{SUVR-FK mapping procedure}
\label{sec:mapping}
In this section, we introduce a validation strategy to map and compare PET images data from \cite{therriault_biomarker_2022} and the numerical results. In \cite{therriault_biomarker_2022},
Therriault et al. conducted a longitudinal \textit{in--vivo} study involving a well--characterized cohort of participants, including cognitively unimpaired individuals (CU), patients with mild cognitive impairment (MCI), and those with AD. Each participant underwent baseline PET scans using [\textsuperscript{18}F]AZD4694 for A\textbeta{,} with SUVR values computed using the cerebellar crus I gray matter as the reference region. Additionally, a PET scan with [\textsuperscript{18}F]MK6240 for \tauP was performed, where SUVR was calculated using the whole cerebellar gray matter as reference. The results were summarized by plotting mean PET-SUVR values across anatomical regions as a function of Braak stages.  
\par
\begin{figure}[t]
    \centering
    \includegraphics[width=\textwidth]{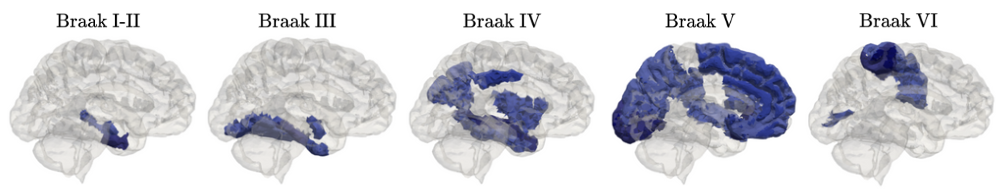}
    \caption{Anatomical representation of the Braak stages used in this work for the left hemisphere of the brain. Each subfigure highlights the brain regions for the specific Braak stage (blue), overlaid on a three--dimensional geometry of the left hemisphere for anatomical reference.}
    \label{fig:Braak_anatomia}
\end{figure}

\begin{table}[t]
\centering
\begin{tabular}{|c|l|}
\hline
\textbf{Stage} & \textbf{Main functions of affected areas} \\
\hline
\textbf{I–II} & Impairments in episodic memory, emotion regulation, and spatial orientation  \\
\hline
\textbf{III–IV} & Decline in language abilities, semantic processing, and attentional control  \\
\hline
\textbf{V–VI} & Loss of executive functions, visuospatial abilities, and visual processing  \\
\hline
\end{tabular}
\caption{Functional interpretation of Braak stages.}
\label{tab:Braak_stages}
\end{table}
To ensure consistency with the methodology used in \cite{therriault_biomarker_2022}, we adopt the same PET-based Braak stage classification in our validation analysis. The stages are qualitatively summarized in Table~\ref{tab:Braak_stages} in terms of the main affected functions. Braak stages were defined according to specific anatomical regions, grouped as follows:
\begin{itemize}
    \item \textbf{Braak I}: transentorhinal cortex;
    \item \textbf{Braak II}: entorhinal cortex and hippocampus;
    \item \textbf{Braak III}: amygdala, parahippocampal gyrus, fusiform gyrus, and lingual gyrus;
    \item \textbf{Braak IV}: insula, inferior temporal, lateral temporal, posterior cingulate and inferior parietal cortices;
    \item \textbf{Braak V}: orbitofrontal, superior temporal, inferior frontal, cuneus, anterior cingulate, supramarginal gyrus, lateral occipital, precuneus, superior frontal, and rostromedial frontal cortices;
    \item \textbf{Braak VI}: paracentral, postcentral, precentral, and pericalcarine cortices.
\end{itemize}
We report in Figure~\ref{fig:Braak_anatomia}, the corresponding anatomical representation of the stages, highlighting in blue the specific regions affected at each stage. The three-dimensional brain mesh was segmented in Freesurfer \cite{fischl_freesurfer_2012}, using the Desikan-Killiany cortical atlas \cite{desikan_automated_2006}. Because the available parcellation does not include the transentorhinal cortex, the validation analysis is restricted to Braak stages II-VI.
\par
The misfolding proteins that are investigated exhibit distinct spatio-temporal evolution: \Abeta remains primarily confined to the neocortex, whereas \tauP progressively invades different brain regions following specific spatial patterns \cite{whittington_spatiotemporal_2018}.
Accordingly, different validation strategies are adopted for the two proteins, both based on identifying Braak macrostages over the simulation time interval and defining a normalisation function for clinical data.
\subsection*{Mapping procedure for \Abeta}
\label{sec:mapping_abeta}
For \Abeta, simulated concentrations do not permit the establishment of a temporal ordering of Braak macrostages. In the absence of longitudinal studies quantifying transition times between Braak stages with respect to \Abeta accumulation, we proceed as follows. Temporal durations of Braak macrostages are anchored to two landmarks from \cite{therriault_biomarker_2022}: \Abeta positivity (cutoff = 1.55, as average on the neocortex, see \cite{therriault_determining_2021}) from Braak stage II onward, and the onset of the saturation phase from Braak stage IV onward. The PET-SUVR original data $s$ are mapped to the normalized variable $\hat{s}\in[0,1]$ through the following function:
\begin{equation}
\hat{s} =
\begin{cases}
0, & s \leq \theta_{\mathrm{low}}, \\
\dfrac{s - \theta_{\text{low}}}{\theta_{\text{high}} - \theta_{\text{low}}}, & \theta_{\text{low}} < s < \theta_{\text{high}}, \\
1, & s \geq \theta_{\text{high}},
\end{cases}    
\end{equation}
where we fix the lower and upper mapping thresholds $\theta_{\text{low}} = 1.3$ and  $\theta_{\text{high}} = 2.2$, respectively. Those values are computed as the minimum and maximum mean neocortical SUVR values reported in \cite{therriault_biomarker_2022}. This choice allows for an accurate identification of the saturation phase, ensuring that it is neither excessively delayed nor inconsistent with the expected steepness of the sigmoidal curve.
\par
The linear mapping adopted in the activation phase of the protein is consistent with the observation that the \Abeta evolution is well described by a logistic function, due to the minor role of diffusion mechanism in its progression \cite{whittington_spatiotemporal_2018}.
Indeed, since the analysis is based on mean PET-SUVR values and the temporal evolution is subdivided into three macrostages, the linear mapping serves as a suitable approximation of the logistic trend near its inflection point.
\subsection*{Mapping procedure for \tauP}
\label{sec:tau_mapping}
Given the sequential pattern of \tauP propagation across brain regions during disease progression, the Braak stages can be reconstructed directly from the model along the simulated time interval. We observe that the time of maximum production of misfolded protein (inflexion point) for the space average concentration $\langle c_h \rangle$ within each Braak region of interest (ROI) occurs after the activation of the corresponding subregion. Therefore, we propose a critical threshold concentration ($\ccrit$) calibrated using numerical results and clinical data. Moreover, three macrostages are identified along the simulation timeline: Braak Stages 0-II, Braak Stages III-IV, and Braak Stages V-VI. \tauP PET data were reported as the mean SUVR for each Braak region (I to VI), together with region-wise abnormality thresholds defined as 2.5 standard deviations higher than the mean SUVR of CU young adults \cite{pascoal_longitudinal_2021}. 
\par
To normalise the PET-SUVR data, we rely on tau aggregation kinetics, which are characterised by three principal phases: lag, growth, and plateau \cite{kamath_kinetics_2021}. The function is defined as follows:
\begin{equation*}
\hat{s} =
\begin{cases}
0, & \text{Stationary phase }(s \leq \theta_{\mathrm{low}}), \\
\gamma \  \dfrac{s - \theta_{\text{low}}}{\theta_{\text{high}} - \theta_{\text{low}}}, & \text{Lag phase}, \\
\dfrac{s - \theta_{\text{low}}}{\theta_{\text{high}} - \theta_{\text{low}}}, & \text{Active phase}, \\
1, & \text{Saturation phase }(s \geq \theta_{\text{high}}),
\end{cases}
\end{equation*}
where $\gamma\in(0,1)$ is the scaling factor of the lag phase, while $\theta_{\mathrm{low}} = 0.75$ and $\theta_{\mathrm{high}} = 2.20$ are the lower and upper mapping thresholds, respectively. $\theta_{\mathrm{low}}$ is the minimum global mean tau PET--SUVR value, while the upper threshold $\theta_{\text{high}}$ is set as the smallest of the regional maximum values. This choice ensures that all Braak ROIs can reach the saturation phase even after applying the mapping procedure.
\par
The definition of \tauP accumulation phases in each Braak subregion relies on the biphasic SUVR increase pattern observed in Figure 1b \cite{therriault_biomarker_2022}. The active phase begins at the Braak stage immediately preceding the region's crossing its \tauP abnormality threshold. However, if the mean SUVR at that stage falls below this threshold by an empirical value $\epsilon$ (calibrated from clinical data in \cite{therriault_biomarker_2022}), the active phase onset shifts back two Braak stages to account for potential biological initiation during the prior transition. Preceding stages are then classified as \textit{lag} or \textit{stationary phases}. An exception applies to Braak II, an early-affected region, where classification starts directly with the active phase.
\par
The scaling factor $\gamma$ modulates the mapping from the \textit{lag} to \textit{active phase}, since it reflects the slower dynamics in this regime. A value of $\gamma = \frac{5}{20}$ provides a good fit to the PET-SUVR data and is therefore adopted in the following validation results.
\par
The use of linear mappings is supported by \cite{chaggar_personalised_2025}, where a mechanism-based model calibrated on PET data introduces a normalized variable $q_i = (s_i - s_{0,i}) / (s_{\infty,i} - s_{0,i}), \; i = 1,\ldots,R$ to capture regional variability in \tauP aggregation. Here, $s_{0,i}$ represents the SUVR value corresponding to a healthy state in region $i$, while $s_{\infty,i}$ denotes the limiting SUVR value associated with an advanced stage of the pathology. 
That study demonstrates accurate fitting on both in-sample and out-of-sample SUVR data.
Since our analysis is based on mean clinical data, region-specific threshold tuning would not provide significant improvements. Therefore, we apply an analogous mapping with global parameters.
\section{Results}
\label{sec:results}
In this section, we present results obtained with the proposed methods applied to the spread of \Abeta and \tauP in the brain.
In the analysis, we will report the space-averaged protein concentrations. For a generic function $g$ over the continuous brain domain $\Omega$, the space average is defined as: $\langle g \rangle = \frac{1}{|\Omega|}\int_{\Omega} g(\boldsymbol{x})\,\mathrm{d}\boldsymbol{x}$.
In the reduced-order model, the brain domain $\Omega$ is subdivided into a finite partition $\{\Omega_i\}_{i=1}^{M}$, where each subdomain $\Omega_i$ is a node of the cerebral graph and is characterised by a measure $| \Omega_i |$. The corresponding discrete approximation of the space average reads: $ \langle g \rangle=\frac{1}{|\Omega|}\sum_{i=1}^{M} g_i\,|\Omega_i|$.
%
\subsection{Three-dimensional high-fidelity simulations}
In this section, we show the results of the 3D simulations obtained using the high-fidelity discretisation described in Section~\ref{sec:polyDG}.
\par
\begin{figure}[t]
    \centering
    \begin{subfigure}[b]{0.54\linewidth}
        \centering
        \includegraphics[width=\linewidth]{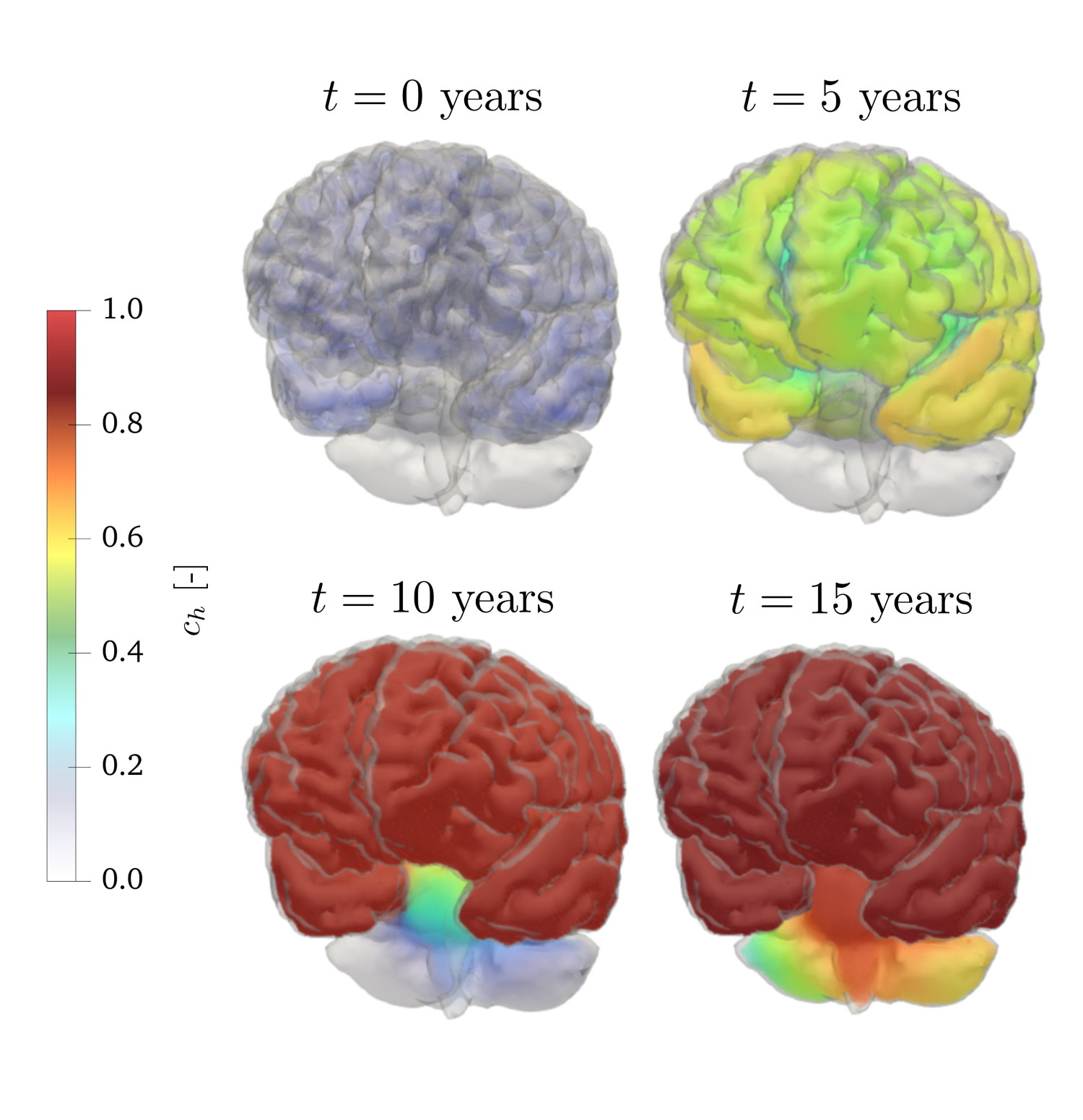}
        \caption{}
        \label{fig:3D_abeta_solution}
    \end{subfigure}
    \begin{subfigure}[b]{0.45\textwidth}
        \centering
        \includegraphics[width=\linewidth]{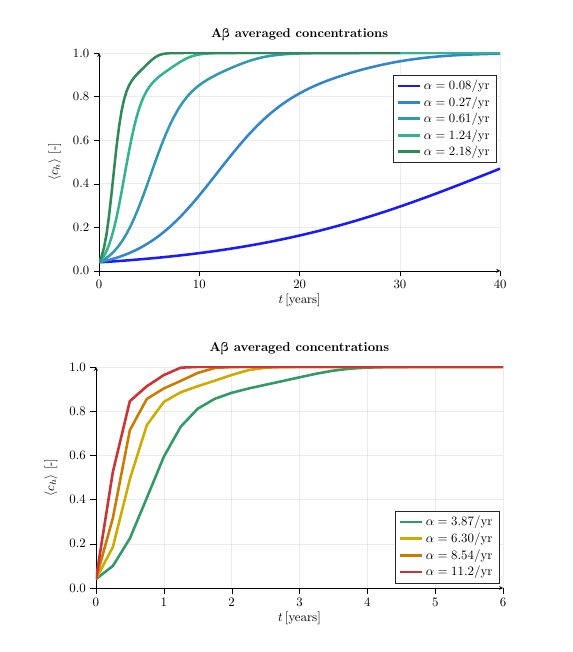}
        \caption{}
        \label{fig:3D_abeta_sensitivity}
    \end{subfigure}
    \caption{\Abeta numerical solution at different stages of the pathology for $\alpha=0.61\,\mathrm{years}^{-1}$ \textbf{(a)} and \Abeta average concentration $\langle c_h\rangle$ w.r.t. time $t$ for different choice of $\alpha$ \textbf{(b)}.}
    \label{fig:3D_abeta_complete}
\end{figure}
Figure~\ref{fig:3D_abeta_solution} shows the simulated spatio-temporal evolution of the \Abeta protein within the brain, obtained by fixing $\alpha = 0.61\,\mathrm{years}^{-1}$, which is the average value obtained by biological measurements in \cite{corti_exploring_2023}. The initial condition collocates the misfolded \Abeta in all the upper cortex (see time $t=0$ years). Moreover, using all the values adopted in \cite{corti_exploring_2023}, we extend the sensitivity analysis to the 3D geometry. The results are summarised in Figure~\ref{fig:3D_abeta_sensitivity}, which reports the temporal evolution of space average \Abeta protein concentration $\langle c_h\rangle$ as $\alpha$ varies. We can observe a faster increase in the protein concentration for higher values of $\alpha$.
\par
\begin{figure}[t]
    \centering
    \begin{subfigure}[b]{0.54\linewidth}
        \centering
        \includegraphics[width=\linewidth]{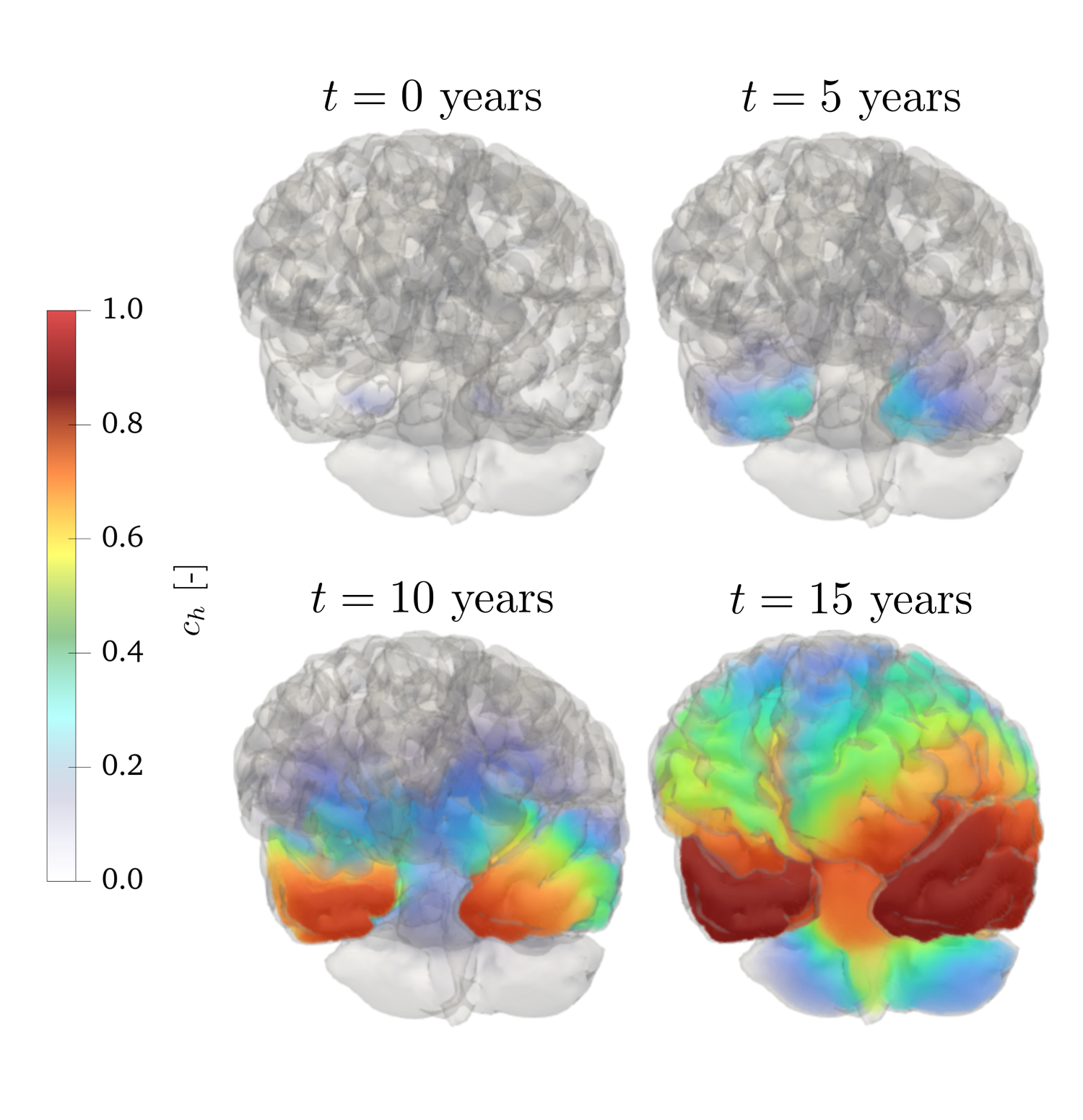}
        \caption{}
        \label{fig:3D_tau_solution}
    \end{subfigure}
    \begin{subfigure}[b]{0.45\textwidth}
        \centering
        \includegraphics[width=\linewidth]{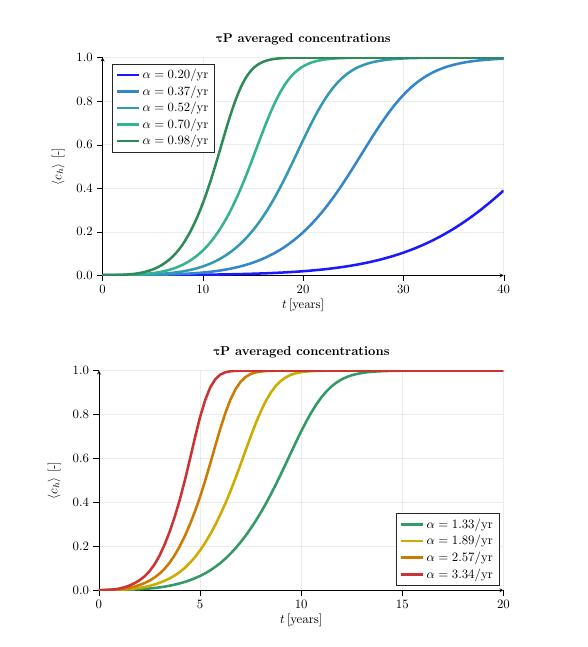}
        \caption{}
        \label{fig:3D_tau_sensitivity}
    \end{subfigure}
    \caption{\tauP numerical solution at different stages of the pathology for $\alpha=0.70\,\mathrm{years}^{-1}$ \textbf{(a)} and \tauP average concentration $\langle c_h\rangle$ w.r.t. time $t$ for different choice of $\alpha$ \textbf{(b)}.}
    \label{fig:3D_tau_complete}
\end{figure}
Analogously, in Figure~\ref{fig:3D_tau_solution}, we show the numerical results obtained for the \tauP protein. In this case, the initial misfolded concentration is located in the enthorinal cortex. The numerical simulation is performed with $\alpha=0.70\,\mathrm{years}^{-1}$. In Figure~\ref{fig:3D_tau_sensitivity}, we report the results of the sensitivity analysis to investigate the influence of the conversion coefficient on the progression of the pathology. We can observe a slower dynamics than in the \Abeta simulations.
\subsection{High-fidelity model validation with PET-SUVR data}
In this section, we validate the high-fidelity simulations using the PET-SUVR data from \cite{therriault_biomarker_2022}. In particular, we consider three different values of the conversion coefficient $\alpha$, using the mean values and the two symmetric ones around it from \cite{corti_discontinuous_2023}. 
\begin{figure}[t]
    \centering
    \includegraphics[width=\linewidth]{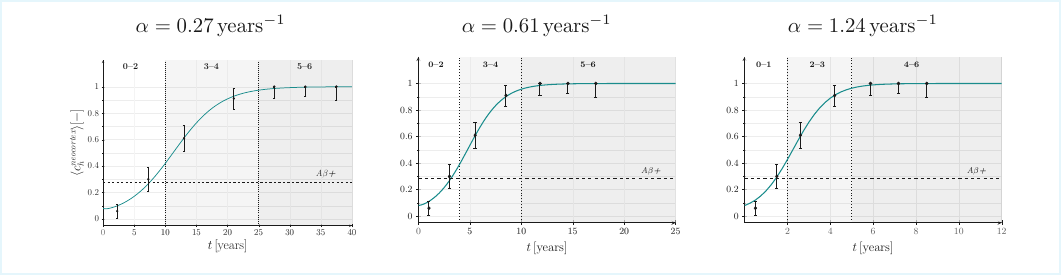}
    \caption{Temporal evolution of the space average \Abeta concentration in the neocortex for three values of $\alpha$. Clinical data with error bars overlaid \cite{therriault_biomarker_2022}.}
    \label{fig:AB_validation}
\end{figure}
\begin{figure}[t]
  \centering
  \begin{subfigure}[b]{0.59\textwidth}
    \centering
    \includegraphics[width=\textwidth]{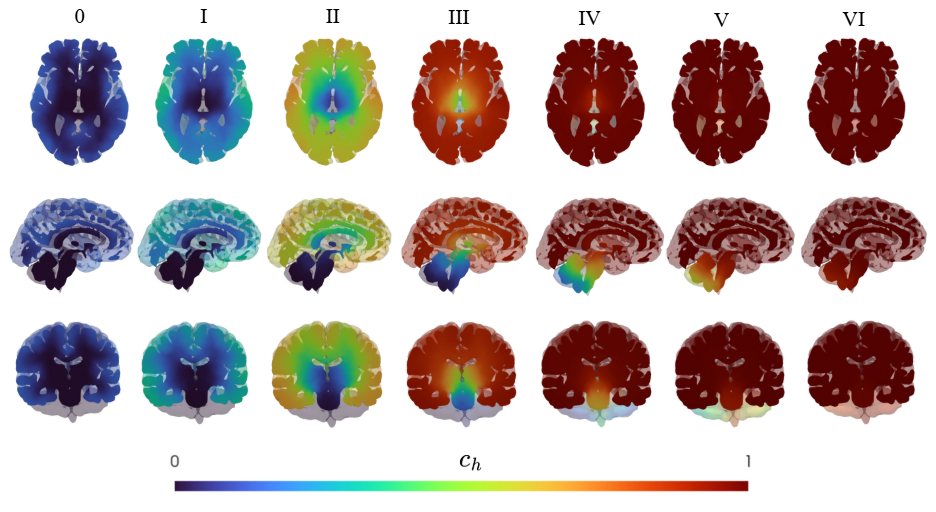}
    \caption{}
    \label{fig:AB061_vs_Bstages_a}
  \end{subfigure}
  \begin{subfigure}[b]{0.39\textwidth}
    \centering
    \resizebox{\textwidth}{!}{\input{Images/AB061_vs_Braakstages}}
    \caption{}
    \label{fig:AB061_vs_Bstages_b}
  \end{subfigure}
  \caption{\textbf{(a)} Spatial distribution of simulated \Abeta concentration at Braak stages I-VI for $\alpha = 0.61\,\text{year}^{-1}$, shown in transversal, sagittal, and coronal views. \textbf{(b)} Simulated spatial average \Abeta concentration across Braak stages.}
  \label{fig:AB061_vs_Bstages}
\end{figure}
\paragraph{Validation of \Abeta simulations.} \Abeta validation results are shown in Figure~\ref{fig:AB_validation}, where simulated curves of average concentration of \Abeta in the neocortex are compared with clinical data from \cite{therriault_biomarker_2022}. The analysis is expressed as a function of Braak macrostages (see Figure~\ref{fig:AB061_vs_Bstages_b}). Moreover, for completeness, Figure~\ref{fig:AB061_vs_Bstages_a} shows an example of the spatial distribution within the simulated brain domain at fixed Braak stages. 
The simulated \Abeta curves are consistent with the trends observed in the clinical data of \cite{therriault_biomarker_2022}. The simulation curve is always within the confidence bands of the data after applying the SUVR mapping described in Section~\ref{sec:mapping_abeta}. As shown in Figure~\ref{fig:AB_validation}, the curves reproduce the clinical data consistently as $\alpha$ varies, except for a temporal shift.
\paragraph{Validation of \tauP simulations.} The validation of \tauP simulations is presented for two representative values of the critical threshold, $\ccrit = 0.5$ and $\ccrit = 0.4$, shown in Figure~\ref{fig:tau_validation_crit05} and Figure~\ref{fig:tau_validation_crit04}, respectively. In both figures, each row corresponds to a Braak region, and the columns show the validation outcomes for different values of the conversion coefficient $\alpha$.
\begin{figure}[ht!]
\centering
\includegraphics[width=\linewidth]{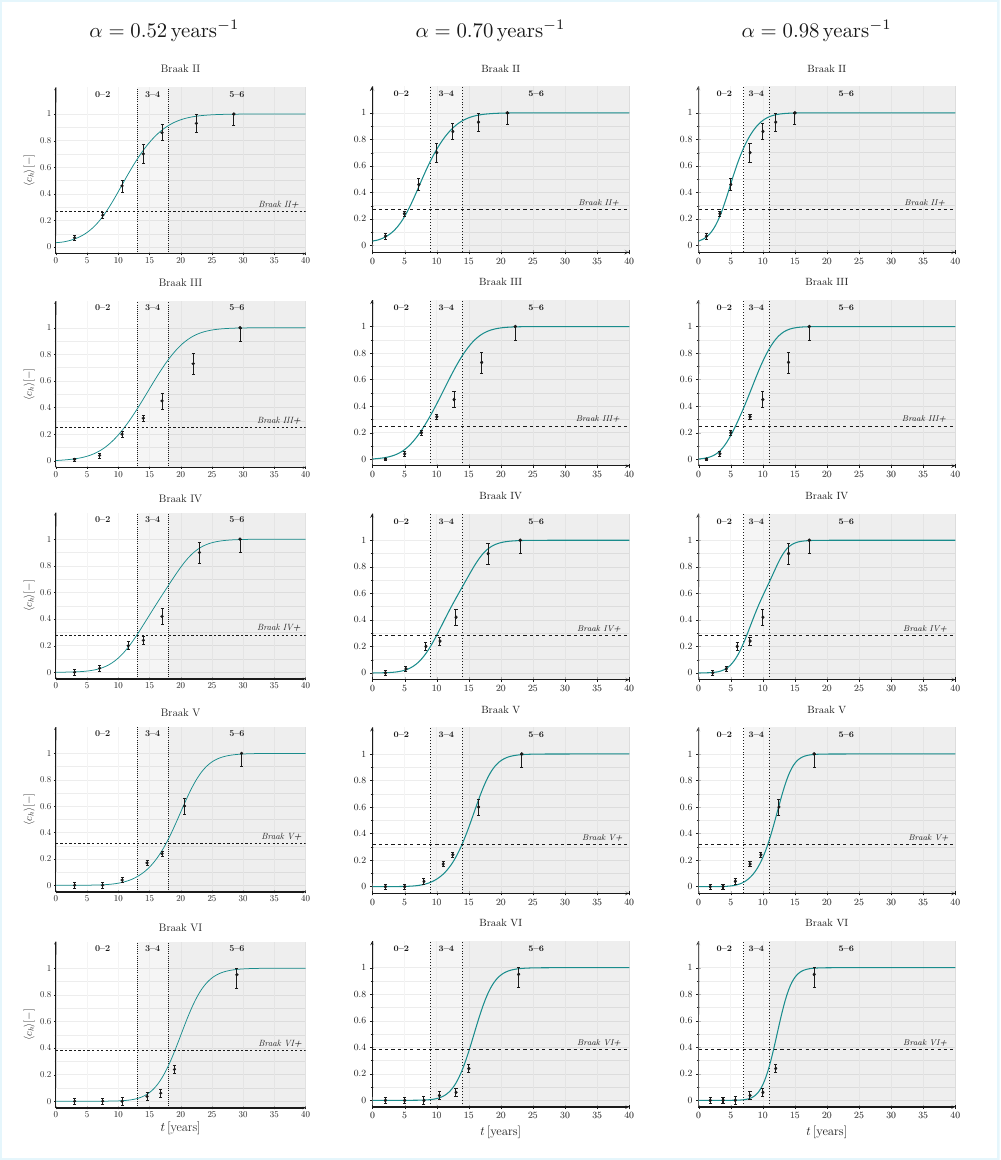}
\caption{Validation results. The solid lines represent the numerically simulated spatial averages of \tauP concentration in Braak ROIs, while the black dots with error bars indicate the clinical PET-SUVR data \cite{therriault_biomarker_2022}. The critical threshold $\ccrit=0.5$ is used to reconstruct the Braak macrostages over the simulated time interval.}
\label{fig:tau_validation_crit05}
\end{figure}
\begin{figure}[ht!]
\centering
\includegraphics[width=\linewidth]{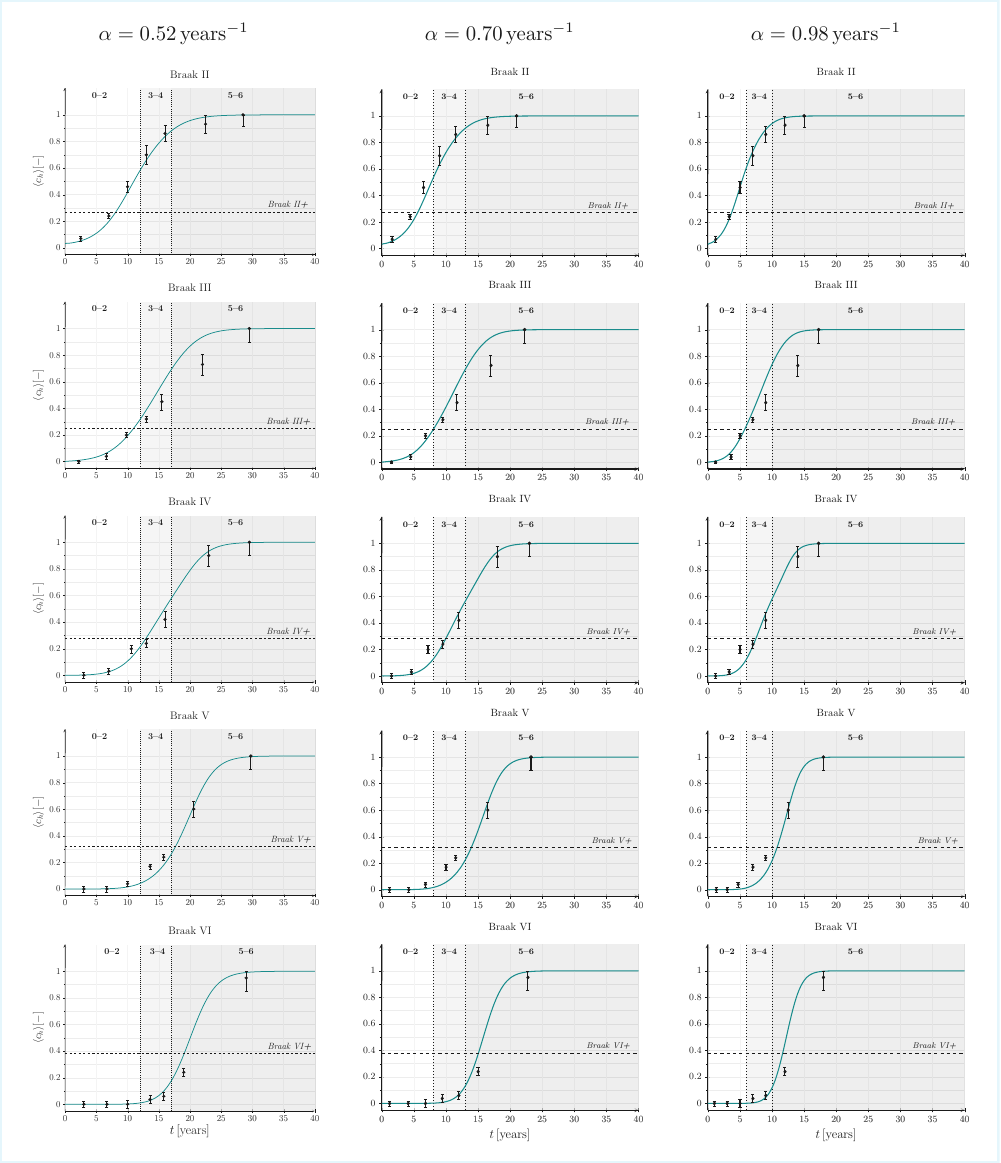}
\caption{Validation results. The solid lines represent the numerically simulated spatial averages of \tauP concentration in Braak ROIs, while the black dots with error bars indicate the clinical PET-SUVR data \cite{therriault_biomarker_2022}. The critical threshold $\ccrit=0.4$ is used to reconstruct the Braak macrostages over the simulated time interval.}
\label{fig:tau_validation_crit04}
\end{figure}
The numerical simulations of \tauP exhibit good agreement with the clinical data reported in~\cite{therriault_biomarker_2022} for both critical threshold values. In Figure~\ref{fig:tau_validation_crit05}, a slight underfitting is observed in Braak~III during the transition between the last two macrostages. When the threshold is reduced to $\ccrit = 0.4$ (Figure~\ref{fig:tau_validation_crit04}), this discrepancy is reduced. For this latter case, the agreement with the clinical data is satisfactory across all curves.
The critical threshold $\ccrit$ affects the temporal positioning of the Braak macrostages in post-processing. The conversion coefficient $\alpha$ governs the rate of disease progression. The duration of the central macrostage remains approximately constant across both threshold values, although temporally shifted.
\subsection{Comparison of high-fidelity and reduced-order descriptions}
In this section, we report the results of the comparison between the high-fidelity 3D model, the corresponding 2D simulations on the sagittal section from \cite{corti_exploring_2023}, and the reduced-order graph-based model. This comparison is fundamental to evaluating how the computational efficiency gains resulting from the geometrical reduction affect the description of protein spreading in AD.
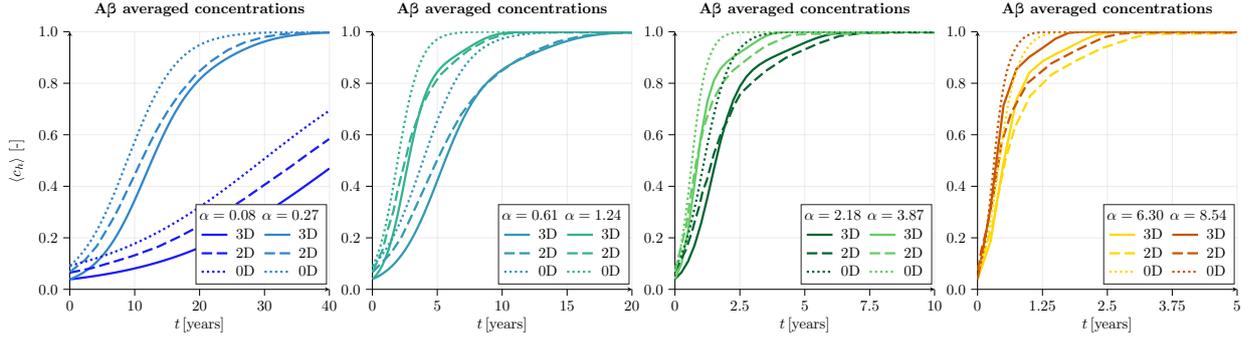
\begin{figure}[t]
    \resizebox{\textwidth}{!}{\input{Tikz/AB_3D2Dgraph}}
    \caption{Computed \Abeta average concentration $\langle c_h\rangle$ as a function of time $t$ for different choices of $\alpha\;[\mathrm{years}^{-1}]$, comparing high-fidelity (3D and 2D) and low fidelity results.}
    \label{fig:AB_3D2Dgraph}
\end{figure}
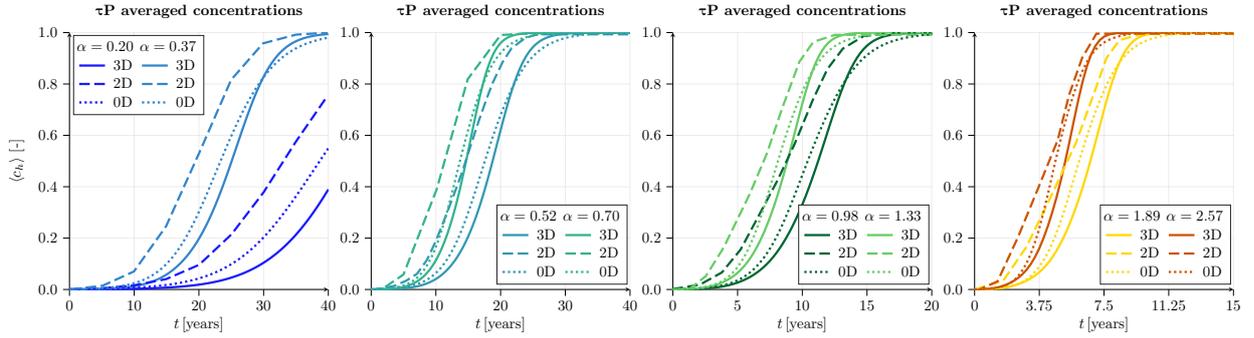
\begin{figure}[t]
    \resizebox{\textwidth}{!}{\input{Tikz/Tau_3D2Dgraph}}
    \caption{Computed \tauP average concentration $\langle c_h\rangle$ as a function of time $t$ for different choices of $\alpha\;[\mathrm{years}^{-1}]$, comparing high-fidelity (3D and 2D) and low fidelity results.}
    \label{fig:Tau_3D2Dgraph}
\end{figure}
\par
Figures~\ref{fig:AB_3D2Dgraph} and~\ref{fig:Tau_3D2Dgraph} show the spatially averaged concentrations of \Abeta and \tauP over time, respectively, for different values of the parameter $\alpha$ used in the sensitivity analysis. The comparison among different geometrical representations reveals distinct differences in protein concentration dynamics. In particular, the simulations in the 2D sagittal section anticipate the temporal evolution of both \Abeta and \tauP, with the difference more pronounced for \tauP than for \Abeta. Notably, the \Abeta 2D simulations exhibit slower dynamics than the 3D model at high values of $\alpha$. The graph-based model exhibits two distinct behaviours for the two proteins. It anticipates \Abeta spreading relative to both 3D and 2D models; conversely, the temporal evolution of \tauP in the graph model is more consistent with the high-fidelity 3D simulations than with the 2D simulations.
\subsection{Comparison in terms of Braak staging theory}
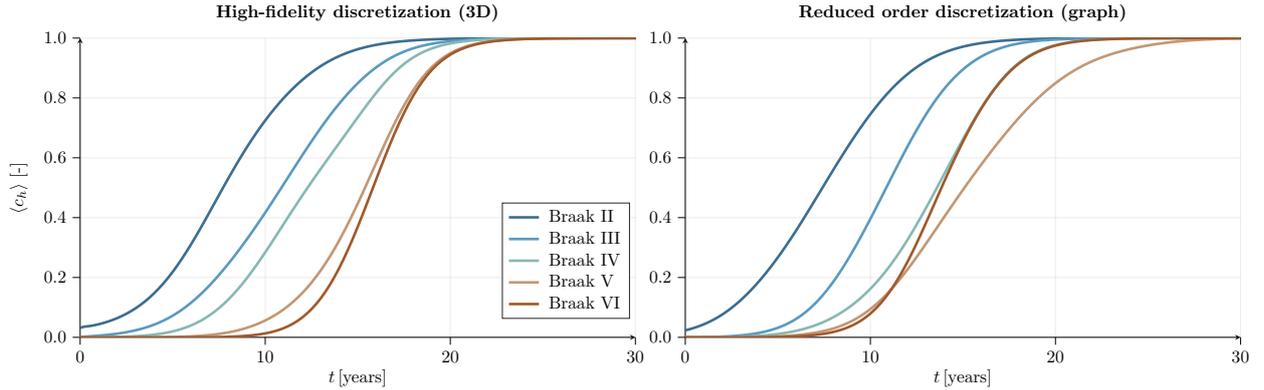
\begin{figure}[t]
    \centering
    \resizebox{\textwidth}{!}{\input{Tikz/BraakStages_tau_Comparison}}
    \caption{Temporal evolution of spatial average of \tauP concentration in the Braak ROIs, using the high-fidelity model (left) and the reduced-order model (right) with $\alpha = 0.70\,\text{year}^{-1}$}
    \label{fig:Braakstages}
\end{figure}
To validate the numerical methods and assess the limitations of graph-based representations of \tauP propagation, we reconstruct the Braak stages by averaging the high-fidelity solver solutions and comparing them with the reduced-order graph-based model. The analysis employs a one-to-one correspondence between three-dimensional Braak regions and graph nodes to ensure consistency in the spatial integration of tau dynamics. Results are presented in Figure~\ref{fig:Braakstages} for $\alpha = 0.70\,\text{year}^{-1}$.
\par
The graph-based model captures the main spatio-temporal features of \tauP propagation through Braak stage IV. However, it does not fully reproduce the progressive invasion of specific anatomical brain regions that characterises AD progression beyond Braak stage IV. By contrast, the correct progression is captured by the high-fidelity 3D simulation \cite{braak_staging_2006}.
\section{Discussion}
In this section, we provide a detailed discussion of the results obtained in the previous section.
\label{sec:discussion}
\subsection{3D high-fidelity results}
The 3D results illustrate the distinct spatio-temporal dynamics of the two proteins. Figure~\ref{fig:3D_abeta_solution} shows that \Abeta rapidly progresses towards the pathological state, in all the brain geometries, mostly due to the diffused initial condition inside the cortex. Concerning \tauP, starting from a localised condition, it spreads through the brain following a specific spatial pattern consistent with the Braak staging theory \cite{whittington_spatiotemporal_2018}. Indeed, Figure~\ref{fig:3D_tau_solution} shows an advanced involvement of the temporal lobe after $10$ years.
Due to the synthetic nature of the FK model, the conversion coefficient $\alpha$ cannot be directly measured from clinical tests. However, its value is of primary importance because it affects the rate of disease evolution, which varies across individuals \cite{fornari_prion-like_2019,corti_exploring_2023}. The sensitivity analysis based on the biological measured values as in \cite{corti_discontinuous_2023} shows that a characteristic sigmoidal behaviour can be observed for both proteins and for any value of $\alpha$. An increase in the conversion coefficient $\alpha$ leads to a sharper initial front and a faster transition to the saturation phase, indicating a more aggressive pathological progression, consistent with previous studies on graphs \cite{fornari_prion-like_2019} or 2D geometries \cite{corti_exploring_2023}. Although this trend is common to both proteins, the dynamics appear to be slower for the \tauP protein, due to the localised initial condition in the brain geometry.
\subsection{High-fidelity model validation with PET-SUVR data}
\paragraph{Validation of \Abeta simulations} The robustness of the model validation across different values of $\alpha$ indicates that the simulated \Abeta kinetics are clinically accurate independently of the temporal scaling of the disease progression. Since the ability of the model to reproduce clinical data is independent of $\alpha$, all simulations remain within the confidence bands of the clinical measurements when appropriately mapped through the SUVR transformation. This consistency across the parameter space demonstrates the reliability of the high-fidelity model in capturing the underlying pathological mechanisms of \Abeta deposition, regardless of the disease development timeline.
To isolate the clinically validated stage-dependent behaviour of A\textbeta{,} the analysis is restricted to Braak macrostages (Figure~\ref{fig:AB061_vs_Bstages_b}), which are defined through external clinical landmarks (positivity cutoff and plateau phase). By expressing the results in terms of Braak stages rather than absolute time, the temporal hypothesis is relaxed, and the intrinsic relationship between pathology progression and brain staging emerges independently of $\alpha$. From a medical perspective, $\alpha$ encodes information not only about the rate of disease progression, but also about the spatial distribution of misfolded protein within the brain (see Figure~\ref{fig:AB061_vs_Bstages_a}). Therefore, Figures~\ref{fig:AB_validation} and~\ref{fig:AB061_vs_Bstages_b} should be considered together: the former validates the model against the clinical data in \cite{therriault_biomarker_2022} for each $\alpha$, confirming that the simulated temporal trajectories are clinically consistent, while the latter reveals the stage-based representation that collapses these trajectories onto a universal curve aligned with clinical neuropathological staging \cite{braak_staging_2006}.

\paragraph{Validation of \tauP simulations} The choice of the conversion coefficient and the critical threshold has a large impact on the final result. The conversion coefficient $\alpha$ governs the rate of pathological progression, while the critical threshold $\ccrit$ determines the onset of local regional activation and marks the transition to the activation phase as its value decreases. An appropriate combination of these two parameters must be identified to optimise the model.
The observation that the transition phase spans approximately 4--5~years, independently of the selected threshold, suggests a robust inherent timescale in the disease progression dynamics. This consistency indicates that the model captures a fundamental aspect of tau pathology progression. These findings are consistent with the longitudinal study by~\cite{krishnadas_rates_2023}, who analysed PET-SUVR data acquired using the tracer [\textsuperscript{18}F]MK6240 for tau, employing the cerebellar cortex as the reference region---a choice analogous to that adopted in~\cite{therriault_biomarker_2022}. In that study, the temporoparietal region became compromised approximately six years after the mesial temporal region, whereas the rest of the neocortex became involved approximately seven years later. A qualitative comparison between these regions and the Braak macrostages employed in the present work reveals a coherent temporal alignment. This correspondence supports the validity and biological plausibility of the simulated progression pattern, demonstrating that our model captures the expected spatio-temporal dynamics of tau pathology observed in clinical studies.
\subsection{Comparison of high-fidelity and reduced-order descriptions}
The accelerated dynamics observed in the 2D model relative to the 3D case likely result from neglecting the brain geometry in the direction orthogonal to the sagittal plane, which leads to a loss of diffusive dynamics and effectively reduces the spatial domain. The larger discrepancies observed for \tauP compared to \Abeta can be attributed to the localisation of the initial seeding: \tauP pathology originates in the entorhinal cortex, a spatially confined region whose representation in a single 2D section is particularly challenging (see Figure~\ref{fig:Tau_3D2Dgraph}). In contrast, \Abeta pathology initiates from more distributed cortical regions, making its projection onto a sagittal section more representative of the full 3D geometry (see Figure~\ref{fig:AB_3D2Dgraph}).
\par
The graph-based model shows complementary strengths and limitations. The accelerated \Abeta spreading reflects the combined effect of a spatially distributed initial condition across the connectome and the simplified geometrical description, which together enhance the effective diffusion of pathology. For \tauP, however, the connectome-based representation captures the temporal dynamics more faithfully than the 2D model. This suggests that, while the graph model reduces spatial resolution at the regional level, it compensates by incorporating long-range axonal connectivity patterns within anatomical regions absent from planar representations. The 2D domain, constrained to a single anatomical plane, cannot recover these network effects. Consequently, when combined with a localised initial condition such as that of \tauP, the graph model describes the temporal evolution more consistently with the 3D physics-based simulations.
\subsection{Comparison in terms of Braak staging theory}
The inability to derive Braak regions directly from the 2D sagittal section stems from the inherent limitation of a single cross-section in representing the complex three-dimensional anatomical structure of the brain. The reduced-order graph-based model successfully captures the main propagation features; however, the simplification of spatial dimension leads to limitations in reproducing regional heterogeneity.
\par
A key factor constraining the model's reliability is the assumption of a constant conversion coefficient $\alpha$. This simplification may be more pronounced when the spatial dimensionality is reduced, as regional heterogeneity is neglected~\cite{corti_uncertainty_2024}. Under this assumption, the reduced-order model with a constant conversion rate appears insufficient to provide a fully accurate representation of Braak stage progression in \tauP propagation, particularly for advanced stages.
\par
The model predictions also depend on the parcellation scheme used to define the network domain. The choice of parcellation can significantly influence both the network topology and quantitative graph measures~\cite{vanwijk_comparing_2010,xu_consistency_2023}, thereby affecting the model's predictions and overall biological accuracy. Future refinements should address parcellation sensitivity and the spatially varying nature of tau conversion dynamics to improve the model's capability to represent the complete Braak staging progression.
\subsection{Limitations}
The main limitations of this work concern the modelling assumptions underlying the conversion coefficient. The assumption of a constant conversion coefficient represents a significant simplification that neglects both spatial and temporal heterogeneity of protein misfolding dynamics. This approximation has a limited impact in the 3D setting, where geometric features and diffusion gradients partially compensate for parameter uniformity; however, its effects become more pronounced in geometrically reduced models. A natural extension would involve the calibration of spatially-localised $\alpha$ parameters directly from PET imaging data, following the approach outlined in~\cite{corti_uncertainty_2024}.
\par
The validation strategy for \Abeta relies on the external identification of Braak macrostages based on clinical evidence, rather than on a direct correspondence between simulated protein accumulation and disease staging emerging from the model dynamics. The availability of longitudinal PET data, or the integration of \Abeta accumulation with calibrated \tauP dynamics, would substantially improve the validation framework and enhance confidence in model predictions.
\par
The present work extends the 2D analysis of the FK equation conducted in~\cite{corti_exploring_2023}. However, increasing attention has been devoted in the literature to the heterodimer model due to its superior descriptive capability in representing clearance, conversion, and production processes separately~\cite{fornari_prion-like_2019,antonietti_discontinuous_2024}. Future investigations should explore the applicability of such mechanistic frameworks within the three-dimensional setting employed herein.
\par
Finally, adapting the proposed strategy to PET imaging data acquired with alternative tracers more widely used in clinical practice~\cite{vanoostveen_imaging_2021}, and validating the simulation framework on a large patient-specific imaging database, could substantially advance our understanding of the practical utility of numerical simulations for pathological prediction and disease progression forecasting.
\section{Conclusions}
\label{sec:conclusions}
In this work, we investigated the spreading dynamics of the two proteins associated with Alzheimer's disease, amyloid-\textbeta{} and tau protein, using the Fisher-Kolmogorov model.
A sensitivity analysis with respect to the conversion coefficient $\alpha$ was performed in a three-dimensional brain domain, highlighting its influence on disease progression rates. Then, the three-dimensional results were compared with those obtained from a two-dimensional sagittal brain section and from a reduced graph-based formulation, to assess the impact of geometrical simplifications on model accuracy and computational efficiency. 
The comparison showed that the two-dimensional geometry does not fully capture the global dynamics of both proteins compared to the three-dimensional setting.
A quantitative validation strategy based on PET-SUVR data was developed using [\textsuperscript{18}F]AZD4694 and [\textsuperscript{18}F]MK6240 tracers to assess amyloid-\textbeta{} and tau proteins, respectively.
The results demonstrated that the three-dimensional model successfully reproduces the observed evolution of both proteins. The analysis, repeated for different values of the conversion coefficient, confirmed that variations in this parameter affect the temporal transition between Braak stages and thus the progression rate, while preserving good agreement with clinical data.
The graph-based model anticipates amyloid-\textbeta{} spreading in time, while providing a more accurate description of tau propagation than the two-dimensional configuration, although it does not reproduce the characteristic Braak staging pattern of tau.
\section*{Acknowledgements}
OASIS-3 provided the brain MRI images: Longitudinal Multimodal Neuroimaging: Principal Investigators: T. Benzinger, D. Marcus, J. Morris; NIH P30 AG066444, P50 AG00561, P30 NS09857781, P01 AG026276, P01 AG003991, R01 AG043434, UL1 TR000448, R01 EB009352. AV-45 doses were provided by Avid Radiopharmaceuticals, a wholly-owned subsidiary of Eli Lilly.

\bibliographystyle{hieeetr}
\bibliography{sample}

\end{document}

%% file: Images/AB061_vs_Braakstages.tex
\begin{tikzpicture}
\begin{axis}[
    width= 1\textwidth,
    height=0.8\textwidth,
    xmin=-0.2, xmax=6.3,
    ymin=-0.05, ymax=1.2,
    xtick      = {0,1,2,3,4,5,6},
    ytick      = {0,0.5,1},
    xticklabels= {0,I,II,III,IV,V,VI},
    xticklabel style={font=\scriptsize, yshift=-3pt},
    yticklabel style={font=\scriptsize, xshift=-3pt},
    xlabel     = {\scriptsize PET-based Braak stage},
    ylabel     = {\scriptsize $\langle c_h \rangle [-]$},
    title       = {\Abeta average in the neocortex}
    axis x line=bottom,
    axis y line=left,
    axis line style={black},
    tick align=outside,
    tick style={black},
    grid=none
]

\addplot+[
  only marks,
  color = black,
  mark=*,
  mark options={fill=teal},
  mark size=2.2pt,
  error bars/.cd, y dir=both, y explicit
] table[x=x, y=y,
  y error plus=eyp,   
  y error minus=eym,  
  col sep=space
]{
  x     y       eyp    eym
  0   0.11     0.0   0.
  1   0.27     0.   0.
  2   0.61     0.   0.
  3   0.90     0.   0.
  4   0.985    0.   0.
  5   0.997    0.   0.
  6   0.999    0.   0.
};

\addplot[black, densely dotted, thick, domain=-0.2:6.3] {0.28};

\node[anchor=west] at (axis cs:5.1,0.33)
  {\textcolor{teal}{\tiny \textit{$A\beta^{+}$}}};
\end{axis}
\end{tikzpicture}

%% file: Tikz/AB_3D2Dgraph.tex
\definecolor{cA}{HTML}{1A1AFF} 
\definecolor{cB}{HTML}{3385CC} 
\definecolor{cC}{HTML}{3399B3} 
\definecolor{cD}{HTML}{33B38F} 
\definecolor{cE}{HTML}{006633}
\definecolor{cF}{HTML}{66CC66}
\definecolor{cG}{HTML}{FFD700}
\definecolor{cH}{HTML}{CC5500} 

\pgfplotsset{
  mysolid/.style={no marks, line width=1.5pt, solid},
  mydashed/.style={no marks, line width=1.5pt, dotted,smooth},
  mydashdot/.style={no marks, line width=1.5pt, dash pattern=on 7pt off 3pt},
  /pgf/number format/read comma as period,
  table/trim cells=true
}

\begin{tikzpicture}
\begin{groupplot}[group style={group size=4 by 1}]
\nextgroupplot[
    width       = 3.00in, 
    height      = 3.00in,
    xmin        = 0.00, 
    xmax        = 40.0,
    ymin        = 0.00, 
    ymax        = 1.00,
    ytick       = {0.0, 0.2, 0.4, 0.6, 0.8, 1.0},
    xtick       = {0, 10, 20, 30, 40},
    xlabel      = {$t\,[\mathrm{years}]$},
    ylabel      = {$\langle c_h \rangle$ [-]},
    axis lines  = left,
    tick pos    = left,
    grid        = both,
    grid style  = {opacity=0.3},
    tick align  = outside,
    tick style  = {black},
    title style = {font=\bfseries},
    title       = {\Abeta averaged concentrations},
    legend columns      = 2,
    legend style        = {at={(0.99,0.33)},legend cell align=left, draw=white!15!black, column sep=0.1cm},
    yticklabel style    = {
    /pgf/number format/.cd, 
    fixed, 
    fixed zerofill, 
    precision=1},
]

\addlegendimage{empty legend}
\addlegendentry{\hspace{-0.85cm}$\alpha=0.08$}

\addlegendimage{empty legend}
\addlegendentry{\hspace{-0.85cm}$\alpha=0.27$}

\addplot+[mysolid,  color=cA]
  table[x=Time, y={avg(c)}, col sep=comma]
  {data/media_pesata_abeta_008.csv};
\addlegendentry{3D}

\addplot+[mysolid,  color=cB]
  table[x=Time, y={avg(c)}, col sep=comma]
  {data/media_pesata_abeta_027.csv};
\addlegendentry{3D}

\addplot+[mydashdot, color=cA]
  table[header=false, x index=0, y index=2, col sep=semicolon]
  {data/abeta_2d_008.csv};
\addlegendentry{2D}

\addplot+[mydashdot, color=cB]
  table[header=false, x index=0, y index=2, col sep=semicolon]
  {data/abeta_2d_027.csv};
\addlegendentry{2D}

\addplot+[mydashed,  color=cA]
    table[x=Time_years, y=Weighted_mean, col sep=comma] 
  {data/weighted_abeta_008_Brain.csv};
\addlegendentry{0D}

\addplot+[mydashed,  color=cB]
    table[x=Time_years, y=Weighted_mean, col sep=comma] 
  {data/weighted_abeta_027_Brain.csv};
\addlegendentry{0D}

\nextgroupplot[
    width       = 3.00in, 
    height      = 3.00in,
    xmin        = 0.00, 
    xmax        = 20.0,
    ymin        = 0.00, 
    ymax        = 1.00,
    ytick       = {0.0, 0.2, 0.4, 0.6, 0.8, 1.0},
    xtick       = {0, 5, 10, 15, 20},
    xlabel      = {$t\,[\mathrm{years}]$},
    axis lines  = left,
    tick pos    = left,
    grid        = both,
    grid style  = {opacity=0.3},
    tick align  = outside,
    tick style  = {black},
    title style = {font=\bfseries},
    title       = {\Abeta averaged concentrations},
    legend columns      = 2,
    legend style        = {at={(0.99,0.33)},legend cell align=left, draw=white!15!black, column sep=0.1cm},
    yticklabel style    = {
    /pgf/number format/.cd, 
    fixed, 
    fixed zerofill, 
    precision=1},
]

\addlegendimage{empty legend}
\addlegendentry{\hspace{-0.85cm}$\alpha=0.61$}

\addlegendimage{empty legend}
\addlegendentry{\hspace{-0.85cm}$\alpha=1.24$}

\addplot+[mysolid,  color=cC]
  table[x=Time, y={avg(c)}, col sep=comma]
  {data/media_pesata_abeta_061.csv};
\addlegendentry{3D}

\addplot+[mysolid,  color=cD]
  table[x=Time, y={avg(c)}, col sep=comma]
  {data/media_pesata_abeta_123.csv};
\addlegendentry{3D}

\addplot+[mydashdot, color=cC]
  table[header=false, x index=0, y index=2, col sep=semicolon]
  {data/abeta_2d_061.csv};
\addlegendentry{2D}

\addplot+[mydashdot, color=cD,restrict y to domain=0:1]
  table[header=false, x index=0, y index=2,
  col sep=semicolon,
    ]{data/abeta_2d_124.csv};
\addlegendentry{2D}

\addplot+[mydashed,  color=cC]
    table[x=Time_years, y=Weighted_mean, col sep=comma] 
  {data/weighted_abeta_061_Brain.csv};
\addlegendentry{0D}
 
 \addplot+[mydashed,  color=cD]
    table[x=Time_years, y=Weighted_mean, col sep=comma] 
  {data/weighted_abeta_124_Brain.csv};
\addlegendentry{0D}

\nextgroupplot[
    width       = 3.00in, 
    height      = 3.00in,
    xmin        = 0.00, 
    xmax        = 10.0,
    ymin        = 0.00, 
    ymax        = 1.00,
    ytick       = {0.0, 0.2, 0.4, 0.6, 0.8, 1.0},
    xtick       = {0, 2.5, 5, 7.5, 10},
    xlabel      = {$t\,[\mathrm{years}]$},
    axis lines  = left,
    tick pos    = left,
    grid        = both,
    grid style  = {opacity=0.3},
    tick align  = outside,
    tick style  = {black},
    title style = {font=\bfseries},
    title       = {\Abeta averaged concentrations},
    legend columns      = 2,
    legend style        = {at={(0.99,0.33)},legend cell align=left, draw=white!15!black, column sep=0.1cm},
    yticklabel style    = {
    /pgf/number format/.cd, 
    fixed, 
    fixed zerofill, 
    precision=1},
]

\addlegendimage{empty legend}
\addlegendentry{\hspace{-0.85cm}$\alpha=2.18$}

\addlegendimage{empty legend}
\addlegendentry{\hspace{-0.85cm}$\alpha=3.87$}

\addplot+[mysolid,  color=cE]
  table[x=Time, y={avg(c)}, col sep=comma]
  {data/media_pesata_abeta_218.csv};
\addlegendentry{3D}

\addplot+[mysolid,  color=cF]
  table[x=Time, y={avg(c)}, col sep=comma]
  {data/media_pesata_abeta_387.csv};
\addlegendentry{3D}

\addplot+[mydashdot, color=cE]
  table[header=false, x index=0, y index=2, col sep=semicolon]
  {data/abeta_2d_218.csv};
\addlegendentry{2D}

\addplot+[mydashdot, color=cF]
  table[header=false, x index=0, y index=2, col sep=semicolon]
  {data/abeta_2d_387.csv};
\addlegendentry{2D}

\addplot+[mydashed,  color=cE]
    table[x=Time_years, y=Weighted_mean, col sep=comma] 
  {data/weighted_abeta_218_Brain.csv};
\addlegendentry{0D}

\addplot+[mydashed,  color=cF]
    table[x=Time_years, y=Weighted_mean, col sep=comma] 
  {data/weighted_abeta_387_Brain.csv};
\addlegendentry{0D}

\nextgroupplot[
    width       = 3.00in, 
    height      = 3.00in,
    xmin        = 0.00, 
    xmax        = 5.00,
    ymin        = 0.00, 
    ymax        = 1.00,
    ytick       = {0.0, 0.2, 0.4, 0.6, 0.8, 1.0},
    xtick       = {0, 1.25, 2.5, 3.75, 5},
    xlabel      = {$t\,[\mathrm{years}]$},
    axis lines  = left,
    tick pos    = left,
    grid        = both,
    grid style  = {opacity=0.3},
    tick align  = outside,
    tick style  = {black},
    title style = {font=\bfseries},
    title       = {\Abeta averaged concentrations},
    legend columns      = 2,
    legend style        = {at={(0.99,0.33)},legend cell align=left, draw=white!15!black, column sep=0.1cm},
    yticklabel style    = {
    /pgf/number format/.cd, 
    fixed, 
    fixed zerofill, 
    precision=1},
]

\addlegendimage{empty legend}
\addlegendentry{\hspace{-0.85cm}$\alpha=6.30$}

\addlegendimage{empty legend}
\addlegendentry{\hspace{-0.85cm}$\alpha=8.54$}

\addplot+[mysolid,  color=cG]
  table[x=Time, y={avg(c)}, col sep=comma]
  {data/media_pesata_abeta_63.csv};
\addlegendentry{3D}

\addplot+[mysolid,  color=cH]
  table[x=Time, y={avg(c)}, col sep=comma]
  {data/media_pesata_abeta_85.csv};
\addlegendentry{3D}

\addplot+[mydashdot, color=cG]
  table[header=false, x index=0, y index=2, col sep=semicolon]
  {data/abeta_2d_630.csv};
\addlegendentry{2D}

\addplot+[mydashdot, color=cH]
  table[header=false, x index=0, y index=2, col sep=semicolon]
  {data/abeta_2d_854.csv};
\addlegendentry{2D}

\addplot+[mydashed,  color=cG]
    table[x=Time_years, y=Weighted_mean, col sep=comma] 
  {data/weighted_abeta_63_Brain.csv};
\addlegendentry{0D}

\addplot+[mydashed,  color=cH]
    table[x=Time_years, y=Weighted_mean, col sep=comma] 
  {data/weighted_abeta_85_Brain.csv};
\addlegendentry{0D}

\end{groupplot}
\end{tikzpicture}

%% file: Tikz/Tau_3D2Dgraph.tex
\definecolor{cA}{HTML}{1A1AFF} 
\definecolor{cB}{HTML}{3385CC} 
\definecolor{cC}{HTML}{3399B3} 
\definecolor{cD}{HTML}{33B38F} 
\definecolor{cE}{HTML}{006633}
\definecolor{cF}{HTML}{66CC66}
\definecolor{cG}{HTML}{FFD700}
\definecolor{cH}{HTML}{CC5500} 

\pgfplotsset{
  mysolid/.style={no marks, line width=1.5pt, solid},
  mydashed/.style={no marks, line width=1.5pt, dotted,smooth},
  mydashdot/.style={no marks, line width=1.5pt, dash pattern=on 7pt off 3pt},
  /pgf/number format/read comma as period,
  table/trim cells=true
}

\begin{tikzpicture}
\begin{groupplot}[group style={group size=4 by 1}]
\nextgroupplot[
    width       = 3.00in, 
    height      = 3.00in,
    xmin        = 0.00, 
    xmax        = 40.0,
    ymin        = 0.00, 
    ymax        = 1.00,
    ytick       = {0.0, 0.2, 0.4, 0.6, 0.8, 1.0},
    xtick       = {0, 10, 20, 30, 40},
    xlabel      = {$t\,[\mathrm{years}]$},
    ylabel      = {$\langle c_h \rangle$ [-]},
    axis lines  = left,
    tick pos    = left,
    grid        = both,
    grid style  = {opacity=0.3},
    tick align  = outside,
    tick style  = {black},
    title style = {font=\bfseries},
    title       = {\tauP averaged concentrations},
    legend columns      = 2,
    legend style        = {at={(0.52,0.99)},legend cell align=left, draw=white!15!black, column sep=0.1cm},
    yticklabel style    = {
    /pgf/number format/.cd, 
    fixed, 
    fixed zerofill, 
    precision=1},
]

\addlegendimage{empty legend}
\addlegendentry{\hspace{-0.85cm}$\alpha=0.20$}

\addlegendimage{empty legend}
\addlegendentry{\hspace{-0.85cm}$\alpha=0.37$}

\addplot+[mysolid,  color=cA]
  table[x=Time, y={avg(c)}, col sep=comma]
  {data/media_pesata_tau_020.csv};
\addlegendentry{3D}

\addplot+[mysolid,  color=cB]
  table[x=Time, y={avg(c)}, col sep=comma]
  {data/media_pesata_tau_036.csv};
\addlegendentry{3D}

\addplot+[mydashdot, color=cA]
  table[header=false, x index=0, y index=1, col sep=semicolon]
  {data/Tau2d_02.csv};
\addlegendentry{2D}

\addplot+[mydashdot, color=cB]
  table[header=false, x index=0, y index=1, col sep=semicolon]
  {data/Tau2d_036.csv};
\addlegendentry{2D}

\addplot+[mydashed,  color=cA]
    table[x=Time_years, y=Weighted_mean, col sep=comma] 
  {data/weighted_tau_020_Brain.csv};
\addlegendentry{0D}

\addplot+[mydashed,  color=cB]
    table[x=Time_years, y=Weighted_mean, col sep=comma] 
  {data/weighted_tau_037_Brain.csv};
\addlegendentry{0D}

\nextgroupplot[
    width       = 3.00in, 
    height      = 3.00in,
    xmin        = 0.00, 
    xmax        = 40.0,
    ymin        = 0.00, 
    ymax        = 1.00,
    ytick       = {0.0, 0.2, 0.4, 0.6, 0.8, 1.0},
    xtick       = {0, 10, 20, 30, 40},
    xlabel      = {$t\,[\mathrm{years}]$},
    axis lines  = left,
    tick pos    = left,
    grid        = both,
    grid style  = {opacity=0.3},
    tick align  = outside,
    tick style  = {black},
    title style = {font=\bfseries},
    title       = {\tauP averaged concentrations},
    legend columns      = 2,
    legend style        = {at={(0.99,0.33)},legend cell align=left, draw=white!15!black, column sep=0.1cm},
    yticklabel style    = {
    /pgf/number format/.cd, 
    fixed, 
    fixed zerofill, 
    precision=1},
]

\addlegendimage{empty legend}
\addlegendentry{\hspace{-0.85cm}$\alpha=0.52$}

\addlegendimage{empty legend}
\addlegendentry{\hspace{-0.85cm}$\alpha=0.70$}

\addplot+[mysolid,  color=cC]
  table[x=Time, y={avg(c)}, col sep=comma]
  {data/media_pesata_tau_052.csv};
\addlegendentry{3D}

\addplot+[mysolid,  color=cD]
  table[x=Time, y={avg(c)}, col sep=comma]
  {data/media_pesata_tau_070.csv};
\addlegendentry{3D}

\addplot+[mydashdot, color=cC]
  table[header=false, x index=0, y index=1, col sep=semicolon]
  {data/Tau2d_0518.csv};
\addlegendentry{2D}

\addplot+[mydashdot, color=cD,restrict y to domain=0:1]
  table[header=false, x index=0, y index=1, col sep=semicolon,
    ]{data/Tau2d_07.csv};
\addlegendentry{2D}

\addplot+[mydashed,  color=cC]
    table[x=Time_years, y=Weighted_mean, col sep=comma] 
  {data/weighted_tau_052_Brain.csv};
\addlegendentry{0D}
 
 \addplot+[mydashed,  color=cD]
    table[x=Time_years, y=Weighted_mean, col sep=comma] 
  {data/weighted_tau_070_Brain.csv};
\addlegendentry{0D}

\nextgroupplot[
    width       = 3.00in, 
    height      = 3.00in,
    xmin        = 0.00, 
    xmax        = 20.0,
    ymin        = 0.00, 
    ymax        = 1.00,
    ytick       = {0.0, 0.2, 0.4, 0.6, 0.8, 1.0},
    xtick       = {0, 5, 10, 15, 20},
    xlabel      = {$t\,[\mathrm{years}]$},
    axis lines  = left,
    tick pos    = left,
    grid        = both,
    grid style  = {opacity=0.3},
    tick align  = outside,
    tick style  = {black},
    title style = {font=\bfseries},
    title       = {\tauP averaged concentrations},
    legend columns      = 2,
    legend style        = {at={(0.99,0.33)},legend cell align=left, draw=white!15!black, column sep=0.1cm},
    yticklabel style    = {
    /pgf/number format/.cd, 
    fixed, 
    fixed zerofill, 
    precision=1},
]

\addlegendimage{empty legend}
\addlegendentry{\hspace{-0.85cm}$\alpha=0.98$}

\addlegendimage{empty legend}
\addlegendentry{\hspace{-0.85cm}$\alpha=1.33$}

\addplot+[mysolid,  color=cE]
  table[x=Time, y={avg(c)}, col sep=comma]
  {data/media_pesata_tau_097.csv};
\addlegendentry{3D}

\addplot+[mysolid,  color=cF]
  table[x=Time, y={avg(c)}, col sep=comma]
  {data/media_pesata_tau_132.csv};
\addlegendentry{3D}

\addplot+[mydashdot, color=cE]
  table[header=false, x index=0, y index=1, col sep=semicolon]
  {data/Tau2d_097.csv};
\addlegendentry{2D}

\addplot+[mydashdot, color=cF]
  table[header=false, x index=0, y index=2, col sep=semicolon]
  {data/tau_2d_133.csv};
\addlegendentry{2D}

\addplot+[mydashed,  color=cE]
    table[x=Time_years, y=Weighted_mean, col sep=comma] 
  {data/weighted_tau_098_Brain.csv};
\addlegendentry{0D}

\addplot+[mydashed,  color=cF]
    table[x=Time_years, y=Weighted_mean, col sep=comma] 
  {data/weighted_tau_133_Brain.csv};
\addlegendentry{0D}

\nextgroupplot[
    width       = 3.00in, 
    height      = 3.00in,
    xmin        = 0.00, 
    xmax        = 15.00,
    ymin        = 0.00, 
    ymax        = 1.00,
    ytick       = {0.0, 0.2, 0.4, 0.6, 0.8, 1.0},
    xtick       = {0, 3.75, 7.5, 11.25, 15},
    xlabel      = {$t\,[\mathrm{years}]$},
    axis lines  = left,
    tick pos    = left,
    grid        = both,
    grid style  = {opacity=0.3},
    tick align  = outside,
    tick style  = {black},
    title style = {font=\bfseries},
    title       = {\tauP averaged concentrations},
    legend columns      = 2,
    legend style        = {at={(0.99,0.33)},legend cell align=left, draw=white!15!black, column sep=0.1cm},
    yticklabel style    = {
    /pgf/number format/.cd, 
    fixed, 
    fixed zerofill, 
    precision=1},
]

\addlegendimage{empty legend}
\addlegendentry{\hspace{-0.85cm}$\alpha=1.89$}

\addlegendimage{empty legend}
\addlegendentry{\hspace{-0.85cm}$\alpha=2.57$}

\addplot+[mysolid,  color=cG]
  table[x=Time, y={avg(c)}, col sep=comma]
  {data/media_pesata_tau_189.csv};
\addlegendentry{3D}

\addplot+[mysolid,  color=cH]
  table[x=Time, y={avg(c)}, col sep=comma]
  {data/media_pesata_tau_257.csv};
\addlegendentry{3D}

\addplot+[mydashdot, color=cG]  
    table[header=false, x index=0, y index=1, col sep=semicolon]
    {data/Tau2d_189.csv};
\addlegendentry{2D}

\addplot+[mydashdot, color=cH]
  table[header=false, x index=0, y index=1, col sep=semicolon]
  {data/Tau2d_257.csv};
\addlegendentry{2D}

\addplot+[mydashed,  color=cG]
    table[x=Time_years, y=Weighted_mean, col sep=comma] 
  {data/weighted_tau_189_Brain.csv};
\addlegendentry{0D}

\addplot+[mydashed,  color=cH]
    table[x=Time_years, y=Weighted_mean, col sep=comma] 
  {data/weighted_tau_257_Brain.csv};
\addlegendentry{0D}

\end{groupplot}
\end{tikzpicture}

%% file: Tikz/BraakStages_tau_Comparison.tex
\definecolor{braakII}{HTML}{3c6e91}
\definecolor{braakIII}{HTML}{5a9abf}
\definecolor{braakIV}{HTML}{87b8b4}
\definecolor{braakV}{HTML}{c39774}
\definecolor{braakVI}{HTML}{a05a2c}

\pgfplotsset{
  mysolid/.style={no marks, line width=1.5pt, solid},
  mydashed/.style={no marks, line width=1.5pt, dotted,smooth},
  mydashdot/.style={no marks, line width=1.5pt, dash pattern=on 7pt off 3pt},
  /pgf/number format/read comma as period,
  table/trim cells=true
}

\begin{tikzpicture}
\begin{groupplot}[group style={group size=2 by 1}]
\nextgroupplot[
    width       = 5.00in, 
    height      = 3.00in,
    xmin        = 0.00, 
    xmax        = 30.0,
    ymin        = 0.00, 
    ymax        = 1.00,
    ytick       = {0.0, 0.2, 0.4, 0.6, 0.8, 1.0},
    xtick       = {0, 10, 20, 30},
    xlabel      = {$t\,[\mathrm{years}]$},
    ylabel      = {$\langle c_h \rangle$ [-]},
    axis lines  = left,
    tick pos    = left,
    grid        = both,
    grid style  = {opacity=0.3},
    tick align  = outside,
    tick style  = {black},
    title style = {font=\bfseries},
    title       = {High-fidelity discretization (3D)},
    legend style        = {at={(0.99,0.45)},legend cell align=left, draw=white!15!black, column sep=0.1cm},
    yticklabel style    = {
    /pgf/number format/.cd, 
    fixed, 
    fixed zerofill, 
    precision=1},
]

\addplot+[no marks, line width=1.5pt, color=braakII]
  table[x=Time, y={avg(c)}, col sep=comma] {data/media_pesata_B1_t07.csv};
\addlegendentry{Braak II}

\addplot+[no marks, line width=1.5pt, color=braakIII]
  table[x=Time, y={avg(c)}, col sep=comma] {data/media_pesata_B3_t07.csv};
\addlegendentry{Braak III}

\addplot+[no marks, line width=1.5pt, color=braakIV]
  table[x=Time, y={avg(c)}, col sep=comma] {data/media_pesata_B4_t07.csv}; 
\addlegendentry{Braak IV}

\addplot+[no marks, line width=1.5pt, color=braakV]
  table[x=Time, y={avg(c)}, col sep=comma] {data/media_pesata_B5_t07.csv};
\addlegendentry{Braak V}

\addplot+[no marks, line width=1.5pt, color=braakVI]
  table[x=Time, y={avg(c)}, col sep=comma] {data/media_pesata_B6_t07.csv};
\addlegendentry{Braak VI}

\nextgroupplot[
    width       = 5.00in, 
    height      = 3.00in,
    xmin        = 0.00, 
    xmax        = 30.0,
    ymin        = 0.00, 
    ymax        = 1.00,
    ytick       = {0.0, 0.2, 0.4, 0.6, 0.8, 1.0},
    xtick       = {0, 10, 20, 30},
    xlabel      = {$t\,[\mathrm{years}]$},
    axis lines  = left,
    tick pos    = left,
    grid        = both,
    grid style  = {opacity=0.3},
    tick align  = outside,
    tick style  = {black},
    title style = {font=\bfseries},
    title       = {Reduced order discretization (graph)},
    legend columns      = 2,
    legend style        = {at={(0.99,0.33)},legend cell align=left, draw=white!15!black, column sep=0.1cm},
    yticklabel style    = {
    /pgf/number format/.cd, 
    fixed, 
    fixed zerofill, 
    precision=1},
]

\addplot+[no marks, solid, line width=1.5pt, color={rgb,1:red,0.235;green,0.431;blue,0.569}]
    table [x index=0, y index=1, col sep=comma, header=true, unbounded coords=discard] {data/weighted_tau07_Brain_braak2.csv};
\addplot+[no marks, solid, line width=1.5pt, color={rgb,1:red,0.353;green,0.604;blue,0.749}]
    table [x index=0, y index=1, col sep=comma, header=true, unbounded coords=discard] {data/weighted_tau07_Brain_braak3.csv};
\addplot+[no marks, solid, line width=1.5pt, color={rgb,1:red,0.529;green,0.722;blue,0.706}]
    table [x index=0, y index=1, col sep=comma, header=true, unbounded coords=discard] {data/weighted_tau07_Brain_braak4.csv};
\addplot+[no marks, solid, line width=1.5pt, color={rgb,1:red,0.765;green,0.592;blue,0.455}]
    table [x index=0, y index=1, col sep=comma, header=true, unbounded coords=discard] {data/weighted_tau07_Brain_braak5.csv};
\addplot+[no marks, solid, line width=1.5pt, color={rgb,1:red,0.627;green,0.353;blue,0.173}]
    table [x index=0, y index=1, col sep=comma, header=true, unbounded coords=discard] {data/weighted_tau07_Brain_braak6.csv};

\end{groupplot}
\end{tikzpicture}